\documentclass[11pt,twoside,onecolumn]{article}

\usepackage{authblk} 

\usepackage{graphicx}
\graphicspath{ {Vis/} }
\usepackage[hang, small,labelfont=bf,up,textfont=it,up]{caption}
\usepackage{subcaption,dsfont}

\usepackage[usenames, dvipsnames]{color}

\usepackage{indentfirst} 

\usepackage[sc]{mathpazo} 
\usepackage[T1]{fontenc} 
\usepackage{microtype} 

\usepackage[english]{babel} 

\usepackage[hmarginratio=1:1,top=32mm,columnsep=20pt]{geometry} 
\usepackage{booktabs} 

\usepackage{lettrine} 

\usepackage{enumitem} 
\setlist[itemize]{noitemsep} 

\usepackage{titlesec} 

\usepackage{fancyhdr} 
\usepackage{hyperref} 

\usepackage{mathtools} 

\usepackage{amsmath,calc,amsthm}

\DeclareMathOperator*{\argmin}{argmin}

\usepackage{amsfonts}
\usepackage{amsthm}

\theoremstyle{definition}

\usepackage{multicol}
\usepackage{multirow}
\usepackage{amssymb}
\usepackage{array}
\usepackage{tabu}
\usepackage{algpseudocode}
\usepackage[numbers]{natbib}

\theoremstyle{definition}

\theoremstyle{plain}

\usepackage{xcolor}

\usepackage{makecell}
\usepackage{scrextend} 

\usepackage{algorithm}
\usepackage{nomencl}
\usepackage{etoolbox}
\makenomenclature

\makeatletter
\patchcmd{\thenomenclature}{\section*}{\section}{}{}
\makeatother

\begin{document}

\title{Real-time Dispatching and Relocation\\of Emergency Service Engineers} 

\author[1]{A. Pechina}
\author[1]{D. Usanov\thanks{Corresponding author. Tel.:+31(0)20 592 4168.\\
\textit{E-mail addresses}: pechina.anna@gmail.com (A. Pechina), usanov@cwi.nl (D. Usanov), p.m.van.de.ven@cwi.nl (P.M. van de Ven), r.d.van.der.mei@cwi.nl (R.D. van der Mei).}}
\author[1]{P.M. van de Ven}
\author[1, 2]{R.D. van der Mei}
\affil[1]{Centrum Wiskunde \& Informatica, Science Park 123,
1098 XG Amsterdam, The Netherlands}
\affil[2]{Vrije Universiteit Amsterdam, De Boelelaan 1081a,
1081 HV Amsterdam, The Netherlands }

\date{} 


\maketitle
\newpage
\begin{abstract}
Capital goods such as complex medical equipment, trains and manufacturing machinery are essential to their users' business, and thus have stringent up-time requirements. Responsive maintenance is crucial for meeting these requirements, which in turn relies on the timely availability of both spare parts and service engineers. Spare parts management for maintenance is well-studied in the research literature, but managing the service engineers has received relatively little attention. In this paper, we consider a network of geographically distributed capital goods, maintained by a set of service engineers who can respond quickly to machine breakdowns. 
We are interested in the question which service engineers to dispatch to what breakdowns, and how to relocate these engineers to maintain good coverage. We propose and evaluate a range of scalable dispatching and relocation heuristics inspired by the extensive research literature in the domain of emergency medical services. We compare the proposed heuristics against each other using comprehensive simulation experiments, and benchmark the best combination of dispatching and relocation heuristics against the optimal policy. We find that this heuristic performs close to optimal, while easily scaling to realistic-sized networks, making it suitable for practical applications.

\end{abstract}



\section{Introduction}
Capital goods are both expensive and essential to the business of their users, and frequent unplanned downtime may have significant repercussions. Consider for instance the reputation damage suffered by train companies unable to maintain the train schedule, or the potentially life-threatening consequences of a broken MRI scanner. In order to ensure continuous operation of these capital goods, 
manufacturers of such products typically provide post-sale support, which may include installation, warranties, spare parts supply and maintenance services. Providing good post-sale support is an important revenue source and competitive advantage for manufacturers \cite{murthy2014}. In 2006, after-sales services accounted for 40\% of the profit of a sample of 120 large US manufacturing companies~\cite{koudal2006}.

An essential component of post-sale support is \textit{corrective} maintenance, i.e., repairing machines that have suffered breakdowns. In practice, a certain service level (e.g., the percentage of failures that should be fixed within a certain time threshold) is determined in the service contract between a manufacturer and the user. Failure to meet these service levels may result in penalties for the manufacturer. 

To meet these ever-tightening service level agreements for corrective maintenance, a manufacturer must be able to quickly dispatch the necessary resources to the site of a breakdown. Since the capital goods are geographically dispersed (e.g., located in many different hospitals), this requires a carefully planned network of spare part warehouses and service engineers. Establishing such networks requires striking a delicate balance between maintaining customer satisfaction and achieving low operational costs~\cite{murthy2004}: building a large number of warehouses would for instance guarantee fast response times, but this is very costly. 
Similarly, maintaining many service facilities and employing a large number of service engineers is costly for a manufacturer. Thus creating a cost-effective maintenance network requires both careful planning of the service facilities and warehouses, and managing the scarce resources in an efficient manner.


While the problem of managing a spare parts network has been extensively studied in research literature (see, e.g., \cite{houtum2015} for the overview), managing service engineers has received little attention so far. Moreover, certain service networks (e.g., software systems support) do not require spare parts at all. Because of this, we focus here on the problem of how to best manage service engineers.



Service engineers are typically located at geographically dispersed base stations, from which they can be quickly dispatched to the site of a breakdown. Here we assume that the location of these base stations  and the number of service engineers are given, and we are interested in the question of how to manage the service engineers. This includes the following decisions: (i) which service engineer to dispatch to a breakdown; (ii) should we dispatch one right away or wait for a nearby engineer to become available; (iii) to which base station should we send an engineer who just finished a repair; and (iv) once a service engineer is dispatched to repair a machine, should we relocate idle service engineers to improve the coverage of the region.

\begin{figure}[t]
\centering
\includegraphics[width = 0.65\textwidth]{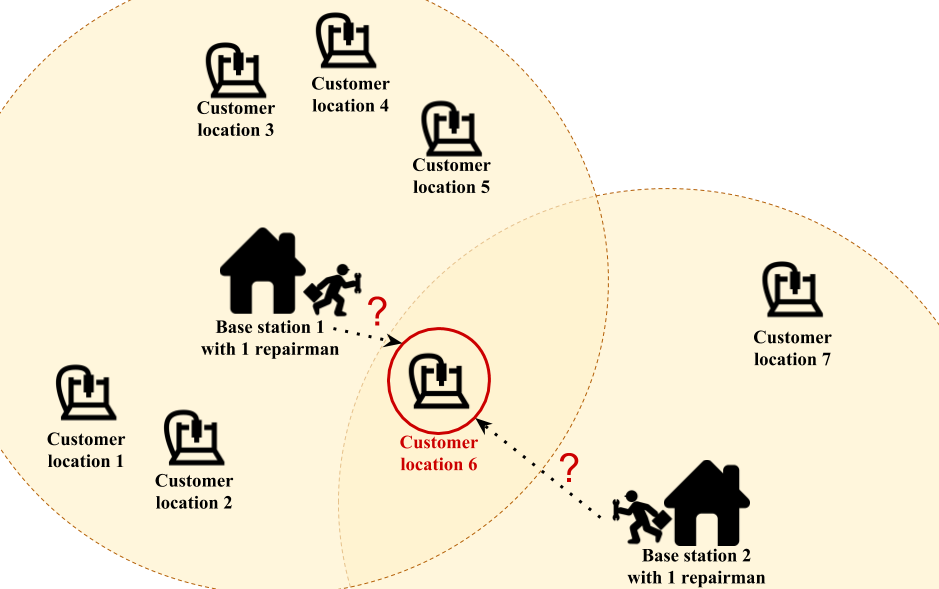}
\captionsetup{width = 0.7\textwidth}
\caption{Example of a dispatching decision}
\label{fig::intoduction}
\end{figure}

Together, these decisions form a trade-off between responding quickly to a current breakdown, and maintaining good coverage for future breakdowns. The goal is to find a dispatching and relocation policy that minimizes the long-term costs caused by violations of the service level agreements with customers. This is a complicated problem due to the fact that dispatching and relocation decisions always increase coverage of one part of the region but decrease coverage of another part. Making a decision requires finding a good balance between the possible costs from failures that are already reported and potential costs from future failures.

Consider, for example, a situation depicted in Figure \ref{fig::intoduction}, where the call arrives from a customer location $6$ and the manager has to decide which service engineer (from base station $1$ or from base station $2$) to dispatch. The first service engineer is closer to the customer and can arrive earlier. However, if the closest service engineer is dispatched, customer locations 1-5 are too far from the remaining idle engineer. If a failure occurs in one of these locations, it can not be fixed in time.

We divide the policy for managing service engineers into two parts: \textit{dispatching} and \textit{relocation}.

\begin{itemize}
\item[--] The dispatching policy is responsible for dispatching service engineers to emergency calls, so it answers the following two questions:
{\it
\begin{itemize}
\item Should a newly reported failure be assigned to one of the idle service engineers or should it be served by one of the busy engineer after he finishes his current job?
\item Which service engineer should be dispatched to the customer?
\end{itemize}
}

\item[--] The relocation policy prescribes the base location of the idle service engineers according to the system state. It answers two questions:
{\it \begin{itemize}
\item Should idle service engineers be relocated from their current base stations to different once to improve coverage?
\item To which base station should we send a service engineer that just finished a job?
\end{itemize}
}
\end{itemize}

Traditionally, Markov decision theory can be used to find the optimal policies. However, for realistic-sized problem instances this approach is infeasible due to computational complexity and high memory usage. Instead, we look for scalable dispatching and relocation heuristics that perform close to optimal. Although management of service engineers is not well-studied in the research literature, we observe that the problem is close to that of Emergency Medical Services (EMS), which also deals with dispatching and relocation of resources~\cite{barneveld2016b}. Recently, much progress has been made in dispatching and relocation for EMS. In this paper we adapt both dispatching and relocation heuristics from that domain to our setting, and show that these perform surprisingly well in this setting as well.

We compare the performance of these heuristics by means of simulation over a wide range of parameters, to identify if there is an approach that performs best for any type of system. In addition, we formulate the model as a Markov decision process, and benchmark the best performing heuristic against the optimal policy for a small instance.

The rest of the paper is organized as follows. Section \ref{sec:lit} gives an overview of related literature. In Section \ref{sec:model}, we present the model and formulate it as a Markov decision process. 
Different dispatching and relocation heuristics are discussed in Sections \ref{sec:disp} and \ref{sec:reloc}, respectively. Numerical experiments are presented in Section \ref{sec:num}. Finally, Section \ref{sec:conclusion} contains conclusions and suggestions for further research.

\section{Literature review}\label{sec:lit}

The problem of real-time management of service engineers arises in the context of spare parts management in a service logistics network. There is an extensive research in the area of spare parts management (see~\cite{houtum2015} for an overview). Most of it focuses on operations for spare parts. One such recent work is by Tiemessen {\em et al.}~\cite{tiemessen2013}, where the authors study the problem of dynamic dispatching of spare parts on a network with multiple customer classes. To our knowledge, there is very limited research on real-time service engineer management. In \cite{benmerzouga1999} and \cite{katehakis1984} the authors study optimal repairmen allocation policies in a simplified setting, for instance ignoring the geographical locations of the engineers. The work closest to ours is by Drent {\em et al.}~\cite{drent2018}, where the authors showed the benefit of deviating from the closest-first dispatching policy, as well as proactively relocating service engineers. However, the authors used a stylized grid-like type of network with deterministic repair and traveling times, where service engineers could reside anywhere on the grid. 

The field of EMS is well-studied and is closely related to our setting. There are however a number of crucial differences between these two application areas. For instance, in our setting there is only a finite number of machines that can break down, and each breakdown affects the rate at which new ones occur. This is in contrast to EMS, where the ongoing incidents do not affect the arrival rate of new incidents. Moreover, many crucial parameters such as the load, coverage area, target response times and service times differ between these two settings, necessitating the present study. It is also worth noting that to our knowledge there is no comprehensive comparison between heuristics in the EMS area, and this study is a first step towards that. In this paper, we will introduce a new model for dynamic management of service engineers, and develop a number of heuristics inspired by the research in the field of EMS. 

In the remainder of this section we outline the most relevant results from the area of EMS operations. For an extensive overview of recent work on location, relocation and dispatching of ambulances, we refer to~\cite{belanger2018}. We organize our literature review according to the type of methods used. First, we discuss Integer Linear Programming (ILP) based approaches. Then we discuss results obtained using Markov Decision Process (MDP) theory, followed by several heuristics that were successfully used for dispatching and relocation. Finally, we cover those references that use Approximate Dynamic Programming (ADP) for dynamic EMS management.

{\bf ILP-based approaches. } A compliance table is a policy that precomputes the optimal locations depending on the number of available service units. Every time this number changes (e.g.,~when a call arrives, or when a service is finished), idle service units are repositioned according to the compliance table. A method called Maximum Expected Coverage Relocation Problem (MEXCRP) was introduced in~\cite{gendreau2006} to compute compliance tables, where for each number of available servers the coverage was maximized. The algorithm was later extended in~\cite{barneveld2016b} to Maximal Expected Coverage Relocation Problem (MEXCRP). MEXCRP compliance tables incorporate the busy fraction of the ambulances.

The problem of choosing a service unit that can be best dispatched to a new emergency call can also be formulated as an ILP problem. According to the computational study of Jagtenberg {\em et al.}~\cite{jagtenberg2017b}, it can even outperform other dispatching policies that use more information about the state of the system.

{\bf MDP-based approaches.} A common way to find the optimal policy is to model the system as an MDP with either continuous or discrete time. For small systems the optimal policy can be found, for instance, using policy or value iteration~\cite{jagtenberg2017a, zhang2008}. For real-life systems, however, the state space is often too large and the problem is computationally intractable. One way to address this problem was considered in~\cite{jagtenberg2017a}: instead of finding the exact optimal policy, the authors perform a limited number of value iteration steps, and compare the results for different numbers of steps with other policies.

{\bf Heuristics.} Another approach to tackle the problem of large state space is to make decisions in real time rather then precomputing the best desicion for each possible state. This applies to both relocation and dispatching problems. The first real-time relocation model was proposed in~\cite{gendreau2001}. It is based on the Double Standard Model~\cite{gendreau1997}, and maximizes the demand covered by at least two vehicles. It also minimizes the relocation costs, so the relocation history is taken into account. Another relocation model, maximizing the preparedness of the system (i.e., the capacity of the system to answer future demands), was introduced in~\cite{andersson2007}. The authors proposed a method to find a relocation policy that minimizes travel times.

In~\cite{jagtenberg2015} the Dynamic Maximal Expected Coverage Location Problem (DMEXCLP) heuristic was proposed for redeployment of service units that just finished their service. The heuristic is based on calculating the expected covered demand and choosing the new location of the service units accordingly. This research was later extended in~\cite{barneveld2016a}, where relocation was allowed not only after the service completion but also right after dispatching a service unit. The authors also studied how different restrictions, such as a restriction on the maximum distance of a relocation, influence the performance of the system. Two types of regions (rural and urban) were considered and it was shown that the optimal strategy depends on the type of the region.

The same ideas can be used for making dispatching decisions. In~\cite{gendreau2001} the authors proposed choosing, among all service unites that can reach the incident location in time, the one that leads to the minimal relocation time. Dispatching a service unit that causes the smallest decrease in preparedness was proposed in ~\cite{andersson2007}. In~\cite{jagtenberg2017b} the expected covered demand was used instead, and the obtained dispatching policy outperformed the commonly used closest-first dispatching policy. Further research also incorporated the possibility of waiting for a busy service unit to finish its service, instead of dispatching an idle one~\cite{barneveld2017}.

{\bf ADP-based approaches.} In~\cite{maxwell2010} an ADP approach using so-called approximate policy iteration was proposed for dynamic ambulance management. The authors considered the problem of ambulance redeployment upon completion of their job, for a system with no other types of relocation and a fixed closest-first dispatching policy. More recently, in~\cite{schmid2012} the same framework was used to optimize both the dispatching policy and redeployment of ambulances upon service completion. Finally, in~\cite{nasrollahzadeh2018} the authors considered a general problem of ambulances dispatching and relocation, with the possibility to reposition idle ambulances and to put incoming calls into a queue.



\section{Model description}\label{sec:model}

In this section we introduce our model, discuss our assumptions, and formulate the problem of finding the optimal relocation and dispatching policy as a continuous-time Markov decision process.

For convenience and ease of interpretation we will refer to the capital goods as machines. We consider a service region that is represented by the set of machines $\mathcal{K}=\{1, \dots, K\}$. Each machine can be either \textit{working}, \textit{in repair}, or \textit{waiting for repair} and we denote by $\kappa_k$ the state of a machine $k \in \mathcal{K}$:
$$
\kappa_k = \begin{cases}
0, \quad \text{if machine $k$ is working;} \\
-1, \quad \text{if machine $k$ is in repair;} \\
t, \quad \text{if machine $k$ has been waiting for repair for time $t$.}
\end{cases}
$$
The state of all machines is described by the vector $\pmb{\kappa} = \left(\kappa_1, \dots, \kappa_K\right)$

All the machines are assumed to be identical, and the time until a working machine breaks down is exponentially distributed with rate $\lambda$. When a new machine breakdown is reported, we need to dispatch exactly one service engineer to its location for repairs. Once the service engineer arrives at the location of the machine, the time to repair a broken machine is exponentially distributed with rate $\mu$. We assume for simplicity that these parameters are the same for all machines, although this assumption can be easily relaxed.

Denote by $\mathcal{M}=\{1, \dots, M\}$ the set of service engineers, and by $\mathcal{R} = \{1,\dots,R\}$ the set of base stations where the service engineers idle when not doing maintenance. The number of service engineers that can stay at the same station is not limited. The locations of all machines and base stations are given, and the traveling time between each pair of locations is known and assumed to be deterministic. The time limit $t^*$ within which a service engineer has to reach a broken machine is given. We assume that the network structure is such that each machine can be reached within the time limit from at least one of the base stations. 

When a new breakdown is reported, it can be either immediately served by one of the idle engineers or it can be put in a queue of pending repair requests to be dealt with later. This is done to allow service engineers that are currently busy but closer to the incident location to respond to the breakdown. A service engineer that is dispatched to a broken machine has to finish a repair before being dispatched or relocated elsewhere. A service engineer traveling to a base station is considered idle and can be dispatched or relocated. A service engineer can be dispatched to a call from the queue in the following two situations: (i) if he/she just finished a repair; or (ii) if he/she is traveling towards a base station.

Relocation can be done either when a service engineer is dispatched to answer a call, or when a repair is finished. This is done to compensate for a possible coverage gap left behind by the dispatched engineer. If a new call is put in the queue, then relocation is not allowed. In the case when a new call is first put in the queue and then one of the traveling repairmen arrived at his destination and is dispatched to this call, the relocation is again not allowed. Only one service engineer can be relocated at a time.


The state of a service engineer $m \in \mathcal{M}$ is represented by the tuple $\mathfrak{m}_m = \left(l_m, d_m\right)$, where $l_m$ denotes the destination (either a machine or a base station) and $d_m$ the time left to reach this destination. If the service engineer is not traveling but is either idle at a base station or repairing a machine, the destination is equal to the current location of the service engineer and the distance is equal to $0$. The state of all service engineers is captured by the vector $\pmb{\mathfrak{m}} = \left(\mathfrak{m}_1, \dots, \mathfrak{m}_M \right)$. 

To describe the process we consider the following event types:
\begin{enumerate}
\item a call arrives from a machine;
\item a repair of a machine is finished;
\item a service engineer arrives at a base station;
\item a service engineer arrives at a machine.
\end{enumerate}

We consider the network state right after one of these events take place, with $e \in \{1,2,3,4\}$ the type of the event. The state of the network can be represented by the tuple $s = \left(t, e, \pmb{\mathfrak{m}}, \pmb{\kappa} \right)$, where $t$ denotes the current time. In the rest of the paper we use $t(s)$, $e(s)$, $l_m(s)$, $d_m(s)$ and $\kappa_k(s)$ to denote the corresponding components of a state $s$. Note that the set $S$ of all possible states of the process is uncountably infinite, and since time is included in the state $s$, it is not recurrent.

\subsection{Action space}\label{sec:actions}
The set of possible actions depends on the state $s$ of the system and is mostly determined by the type of event $e$.  Let $\mathcal{F}\left(s\right) = \{m \in 1,\dots,M \mid l_m(s) \in \mathcal{R}\}$ denote the set of all idle service engineers, and $Q(s)$ the set of all machines in the queue in state $s$. We consider each of the four types of events described in the previous section, and define the corresponding sets of possible actions.\\

\noindent
\textbf{Type 1: A call arrives from machine $k$}

\noindent
In this case, either one of the idle service engineers in $\mathcal{F}\left(s\right)$ can be dispatched to the broken machine, or the call can be placed in the queue. If a service engineer is immediately dispatched, we also allow one idle service engineer to be relocated. Let vector $\pmb{X}$ and matrix $\pmb{Y}$ represent the dispatching and relocation decisions, respectively. Here $X_m = 1$ if service engineer $m \in  \mathcal{F}\left(s\right)$ is assigned to the call, and $X_m = 0$ otherwise. And $Y_{mr}=1$ if service engineer $m \in \mathcal{F}\left({s}\right)$ is relocated to base station $r$, and $Y_{mr}=0$ otherwise. The corresponding action space in state $s$ is given by
\begin{align}
\nonumber \mathcal{A}_1(s) = \bigg\{(\pmb{X}, \pmb{Y}) \mid &\sum_{m \in \mathcal{F}(s)}X_m \le 1; \ \sum_{\substack{m \in \mathcal{F}(s),\\ r \in \mathcal{R}}} Y_{mr}\le 1;  \\ 
&\sum_{\substack{m \in \mathcal{F}(s),\\ r \in \mathcal{R}}} X_mY_{mr}=0; \ \Big(1-\sum_{m \in \mathcal{F}(s)}X_m\Big)\sum_{\substack{m \in \mathcal{F}(s),\\ r \in \mathcal{R}}} Y_{mr} = 0 \bigg\}, \label{eqn:A1}
\end{align}
where the constraints ensure that not more than one service engineer is dispatched, not more than one service engineer is relocated, the relocated service engineer differs from the dispatched one, and no relocation is done if the call is placed in the queue.\\

\noindent
\textbf{Type 2: A repair of machine $k$ is finished by service engineer $m$}

\noindent
Service engineer $m$ can be sent either to one of the base stations, or to one of the machines from the queue. Let $Z_l=1$ represent that the service engineer that just became idle is dispatched to location $l \in Q(s) \cup \mathcal{R}$, and $Z_l=0$ otherwise. Let $Y_{nr}=1$ denote that the service engineer $n \in \mathcal{F}(s)$ is relocated to base station $r$, and $Y_{nr}=0$ otherwise. Then an action is described by the pair $(\pmb{Z}, \pmb{Y})$ of redeployment vector $\pmb{Z}$ and relocation matrix $\pmb{Y}$. The corresponding action space is given by
\begin{equation}\label{eqn:A2}
\mathcal{A}_2(s) = \bigg\{\left(\pmb{Z}, \pmb{Y}\right) \mid  \sum_{l \in Q(s)\cup\mathcal{R}}Z_l = 1; \sum_{\substack{n \in \mathcal{F}\left(s\right),\\ r \in \mathcal{R}}} Y_{nr}\le 1 \bigg\},
\end{equation}
where the constraints ensure, that the service engineer is redeployed to exactly one location, and at most one other service engineer is relocated.\\

\noindent
\textbf{Type 3: Service engineer $m$ arrives at base station $r$}

\noindent
Service engineer $m$ can either be left idle at the station, or dispatched to one of the machines in the queue. Relocations are not allowed in this case. Let $U_k = 1$ if service engineer $m$ is dispatched to machine $k$, and $U_k = 0$ otherwise. Then an action is described by vector $\pmb{U} = (U_1,\dots,U_K)$, and the action space is
\begin{equation}\label{eqn:A3}
\mathcal{A}_3(s) = \bigg\{ \pmb{U} \mid \sum_{k \in Q(s)}U_k \le 1 \bigg\}
\end{equation}
Note that if the queue is empty in state $s$, then so is $\mathcal{A}_3(s)$.

\noindent
\textbf{Type 4: Service engineer $m$ arrives at machine $k$}

\noindent
When the service engineer arrives at machine $k$, the repair process is started, and there are no available actions:
\begin{equation}\label{eqn:A4}
\mathcal{A}_4(s) = \emptyset.
\end{equation}

\subsection{State transitions}\label{sec:trans}

Let $\{s_n\}_{n \in \mathds{N}}$ denote the discrete-time stochastic process of the network state over time, embedded on time instances when an event has been resolved. The evolution of this process is characterized by an action $a_n$, a random element $\omega(s_n, a_n)$, and function $\Phi$ as $s_{n+1} = \Phi\left(s_n, a_n, \omega(s_n, a_n)\right)$. The random element determines the event in state $s_{n+1}$, and depends on both the state and the action. Denote by $d(s_n, a_n)$ the minimum of all non-zero distances remaining for the service engineers to travel after taking action $a_n$ in state $s_n$. If no events of types 1 and 2 occur, then in state $s_{n+1}$ time is equal to $t(s_{n+1}) = t(s_n)+d(s_n, a_n)$, and the event is the arrival of a service engineer with the shortest remaining distance to the destination. If there are no traveling service engineers in the system after taking action $a_n$ in state $s_n$, we set $d(s_n,a_n) = \infty$.

Denote the set of all working machines after taking action $a_n$ in state $s_n$ by $\mathcal{W}(s_n, a_n)$ and the set of all machines in repair by $\mathcal{H}(s_n)$. Let $W(s_n) = |\mathcal{W}(s_n, a_n)|$ and $H(s_n) = |\mathcal{H}(s_n)|$. For each machine in $\mathcal{W}(s_n, a_n)$ the time till breakdown is exponentially distributed with rate $\lambda$. For each machine in $\mathcal{H}(s_n)$ the time till repair is finished is exponentially distributed with rate $\mu$. So the time till the next event of one of the first two types is the minimum of several exponential distributions and is exponentially distributed with rate $\eta(s_n) = \lambda W(s_n)+\mu H(s_n)$. If this time is less than $d(s_n,a_n)$, then the next event is either a call arrival or an end of repair.

Therefore, if in state $s_n$ action $a_n$ is taken, the probability that the next event is
\begin{itemize}
\item the arrival of a call from machine $k \in \mathcal{W}(s_n, a_n)$ is $$\frac{\lambda}{\eta(s_n)}(1-e^{-\eta(s_n) d(s_n,a_n)}),$$ if $\mathcal{W}(s_n, a_n)\ne\emptyset$, and $0$ otherwise;
\item the end of repair of machine $k \in \mathcal{H}(s_n)$ is $$\frac{\mu}{\eta(s_n)}(1-e^{-\eta(s_n) d(s_n,a_n)}),$$ if $\mathcal{H}(s_n)\ne\emptyset$, and $0$ otherwise;
\item the arrival of a service engineer to the destination is $$e^{-\eta(s_n) d(s_n,a_n)},$$ if there are any traveling service engineers, and $0$ otherwise.
\end{itemize}

The time until the next event is distributed as the minimum of an exponentially distributed random variable with rate $\eta(s_n)$ and a constant $d(s_n,a_n)$. When the next event and the time until that event are known, the function $\Phi$ updates the states of all service engineers and machines. The remaining travel times are reduced by the transition time, and the waiting times of the broken machines are increased accordingly. If action $a_n$ includes repositioning of some of the service engineers, then their destinations and remaining distances change according to the action. If the event in state $s_{n+1}$ is the end of repair of machine $k$, then the state of this machine is changed from $-1$ to $0$. If the event in state $s_{n+1}$ is the arrival of a service engineer to machine $k$, then the state of this machine is changed from $\kappa_k(s_n)$ to $-1$.

\subsection{Costs} \label{sec:costs}
In the numerical experiments we compare the performance of different dispatching and relocation policies using the fraction of emergency calls responded to within the time limit $t^*$. This performance measure corresponds to the following cost structure. If a call arrives from machine $k$ and a service engineers does not reach this machine within the time limit $t^*$, then a penalty $1$ is paid. All travel costs and other operational costs are ignored, but could be readily added to the model.

Denote by $c(s_n, a_n, s_{n+1})$ the costs that are charged during the transition period from state $s_n$ to state $s_{n+1}$ when action $a_n$ is taken. The costs are equal to the number of machines who's waiting time exceeds $t^*$ within the interval $\left(t(s_n), t(s_{n+1})\right]$:
$$c(s_n, a_n, s_{n+1}) = \sum_{k=1,\dots,K}\mathbb{I}\{\kappa_k(s_{n+1})\ge t^*\}\mathbb{I}\{\kappa_k(s_{n})<t^*\}.$$

\section{Dispatching heuristics}\label{sec:disp}

The dispatching policy determines {\it when} and {\it which} service engineer to assign to a call, based on the current state of the system. 
When a breakdown occurs, a customer wants a service engineer to arrive on scene as quickly as possible. However, when a service engineer is dispatched from one of the base stations, the coverage of the customers around this base station decreases, which may lead to high future costs. Thus a good dispatching policy finds a balance between minimizing immediate costs and future costs. In order to focus fully on the problem of dispatching, in this section relocation is not allowed. In Section~\ref{sec:reloc} we consider the complementary problem of fixing the dispatching policy and varying the relocation problem. So after the completion of service, the engineer returns to the his/her preallocated base station. The initial allocation of service engineers to base stations is made to maximize the expected covered demand and is a solution to the ILP problem \eqref{eqn::expcovdemalloc} described in Appendix~\ref{apx:opt_allocation}.

In our comparative study we consider the following five dispatching heuristics:
\begin{description}
\setlength{\itemsep}{1pt}
\item[DP1] Closest-first dispatching policy (without waiting);
\item[DP2] Maximal coverage dispatching policy (without waiting);
\item[DP3] Maximal expected coverage dispatching policy (without waiting);
\item[DP4] Minimal response time dispatching policy with unknown remaining repair time;
\item[DP5] Minimal response time dispatching policy with known remaining repair time.
\end{description}

In Section~\ref{sec:disp:nowaiting}, we consider heuristic dispatching policies that put a call into the queue only if there are no idle service engineers. Otherwise, a service engineer has to be dispatched immediately. In Section~\ref{sec:disp:waiting}, we discuss heuristics that allow calls to be put into the queue even if idle engineers are available. A comparison of all dispatching policies based on simulation can be found in Section \ref{sec:num:disp}.

\subsection{Dispatching policies without waiting (DP1, DP2, DP3)} \label{sec:disp:nowaiting}
Recall from~\eqref{eqn:A1} that in the event of type $e(s)=1$, the set of all possible actions is described by a vector $X$ that represents the dispatching decision, and a matrix $Y$ that represents the relocation decision. As in this section relocation is not allowed, all elements of the matrix $Y$ are always set to $0$. We do the same for the relocation matrix $Y$ for $e(s)=2$, see~\eqref{eqn:A2}.

In this section, when a call arrives in the system with at least one idle service engineer, an engineer is dispatched immediately. Recall that $\mathcal{F}(s)$ denotes the set of all idle service engineers in state $s$, so for any state $s$ with $\mathcal{F}(s)\ne \emptyset$ and event $e(s)=1$, the set of possible actions from~\eqref{eqn:A1} reduces to
$$\mathcal{A}_1(s) = \bigg\{(X, Y) \mid  \sum_{m \in \mathcal{F}\left(s\right)}X_m = 1, Y_{mr} =0 \ \ \forall m \in \mathcal{F}\left(s\right), r \in \mathcal{R} \bigg\}.$$
If $\mathcal{F}(s)=\emptyset$, then the call is put in the queue and no decision should be made. Note also that under these restrictions the set of possible actions for $e(s) = 3$ is empty, as a call can not be put in the queue if there is an idle service engineer in the system.

\textbf{DP1}. Consider state $s$ with the event $e(s) = 1$ (i.e., a call arrives from machine $k$) and $\mathcal{F}(s)\ne\emptyset$. To describe the dispatching policy, we need to calculate the vector $X$ depending on state $s$. The simplest and most widely used dispatching policy in practice is the so-called {\bf closest-first} policy. Under this policy, the closest idle service engineer is always dispatched to a call. So
$$X_m = 1 \iff m = \argmin_{n \in \mathcal{F}(s)}\left( ||l_nk||_2 + d_n \right), $$
where $||l_nk||_2$ is the distance in time between the destination of the service engineer $n$ and the source of the call, the machine $k$.


%

\textbf{DP2}. One of the possible metrics of the expected system performance is {\it coverage}, i.e., the number of machines covered by at least one service engineer. When a call arrives, for each idle service engineer $m$ that can reach machine $k$ in time, the remaining coverage of the system without him/her is defined as:
$$coverage(m) = \sum_{k' : \kappa_{k'}=0}\mathbb{I}\left\{\exists n \in \mathcal{F}(s): n\ne m, ||l_mk')||_2 \le t^*\right\}.$$
Then in {\bf coverage-based dispatching}, the service engineer with the biggest remaining coverage is dispatched. If there are no service engineers that can reach the source of the call in time, then the remaining coverage is calculated for all idle service engineers and the one with the biggest remaining coverage is dispatched. In case there are multiple service engineers maximizing the remaining coverage, the closest one is dispatched.


\textbf{DP3}. Note that the coverage only estimates the performance of the system for the next call. Alternatively, one can calculate the {\bf expected covered demand}, the fraction of calls that will be answered in time by the system. As the expected covered demand is hard to compute, for computational study we use the approximation described in Appendix \ref{apx:exp_cov}. Similar to the coverage based approach, first we calculate the remaining expected covered demand after dispatching each of the idle service engineers that can reach the broken machine in time, and dispatch the one with the highest remaining expected covered demand. If there are no idle service engineers that can reach the machine in time, we instead do this procedure for all idle service engineers.


\subsection{Dispatching policies with waiting (DP4, DP5)}\label{sec:disp:waiting}
Under the dispatching policies discussed in the previous section, it is mandatory to dispatch an idle service engineer when a call arrives. However, in practice it may occur that a service engineer close to the new breakdown will finish its repair earlier than any idle service engineer can reach the breakdown, in which case it may be better to wait for the busy service engineer to finish, and dispatch him/her instead.

Inspired by this, we extend the closest-first dispatching policy to the dispatching policy that chooses a service engineer with the smallest response time. For a service engineer $m$ whose destination is a base station $r$ (i.e., $l_m = r$) the response time to a call from machine $k$ is calculated as $rt(k,m) = d_m + ||l_mk||_2$.
For a service engineer $m$ whose destination is machine $k'$ (i.e., $l_m = k'$) the response time consists of the distance left to machine $k'$, the length of repair and the distance from machine $k'$ to machine $k$. Then the response time equals $rt(k,m) = d_m + t_{repair} + ||k'k||_2$. We consider two situations: when the length of repair can be estimated upon arrival of a service engineer to the machine (\textbf{DP5}) and when it stays unknown (\textbf{DP4}). If the length of repair $t_{repair}$ is not known, it can be estimated from its distribution. For the computational study of the second situation it is estimated by an $\alpha^{\text{th}}$ percentile of the repair time distribution. Throughout the computational study we take $\alpha = 80\%$.

In the extended closest-first policy the service engineer $m = \argmin_n rt(k,n)$, that minimizes response time, is assigned to the call. If $m \in \mathcal{F}(s)$ then the service engineer is dispatched immediately. If the service engineer $m$ is busy, then the call is placed in the queue.



\section{Relocation heuristics}\label{sec:reloc}

The relocation policy is responsible for the location of idle service engineers. The simplest relocation policy is the static policy, where each service engineer is assigned to a base station and resides there when idle. In Appendix \ref{apx:opt_allocation} we discuss how to compute the static allocation of the service engineers that maximizes the expected covered demand. We assume that this policy is used when the dispatching policies are studied in isolation. However, when a service engineer is dispatched to a call, a large area of the region may become uncovered and it may be beneficial to reallocate other idle service engineers. In our system we allow one idle service engineer to change the destination either when one of the service engineers is dispatched or finishes a repair.

For our comparative study we choose the following five heuristic relocation policies:
\begin{description}
\item[RP1] Static policy;
\item[RP2] MCRP compliance tables;
\item[RP3] MEXCRP compliance tables;
\item[RP4] DMEXCLP heuristic without constraints;
\item[RP5] DMEXCLP heuristic with constraints.
\end{description}

In this section we consider four relocation policies, in addition to the static policy (RP1). Ppolicies RP2 and RP3 are compliance tables constructed according to two different algorithms. The compliance table relocation policy is the one where location of service engineers depends only on the number of idle service engineers. The locations of the idle engineers are ignored, as well as the state of the machines. This allows us to reduce the number of considered situations and precompute relocation actions for larger systems. We present these approaches in Sections~\ref{sec:reloc:mcrp} and~\ref{sec:reloc:amexprep}, respectively. Policies RP4 and RP5 are heuristic relocation policies based on the DMEXCLP heuristic introduced in~\cite{jagtenberg2015}. The main difference between this heuristic and compliance tables is that decisions are made in real-time, that allows to use the information about the current location of service engineers and the state of the machines. We construct two versions of the DMEXCLP heuristic, with no restrictions and with restrictions on relocation, in Section~\ref{sec:reloc:dmexclp}.

\subsection{MCRP compliance tables (RP2)}\label{sec:reloc:mcrp}

The first type of compliance table we present is the Maximum Coverage Relocation Problem (MCRP) compliance table that aims to maximize the probability that the next call is answered in time. Consider a system with $M$ service engineers. Then a compliance table consists of $M$ levels, where level $m$ contains the allocation solution for $m$ idle service engineers in the system. If one of them is dispatched, other service engineers are relocated according to level $m-1$. If one of service engineer becomes idle after finishing a repair, then the idle service engineers are relocated according to level $m+1$.

Denote by $z_{mk}$ the indicator of the fact that at the level with $m$ idle service engineers the machine $k$ is covered (meaning that at least one service engineer can reach it in time), Then at the level $m$ we want to maximize $\sum_{k=1}^Kz_{mk}$, i.e., the number of covered machines. If $x_{mr}$ is the number of service engineers at the base station $r$ at level $m$, then $z_{mk} \le \sum_{r \in N_k}x_{mr}, \quad k = 1,\dots, K$, where $N_k$ is the set of all base stations from which the machine $k$ can be reached in time.

Recall that we allow to relocate at most one service engineer at a time. To include this restriction in the ILP formulation, we introduce non-negative variables $\alpha_{mr}$ that represent the number of service engineers that arrived at base station $r$ after going from level $m+1$ to level $m$. Then $x_{mr} - x_{m+1,r} \le \alpha_{mr}$, $r = 1, \dots, R,\ m = 1,\dots,M-1$, and $\sum_{r=1}^R \alpha_{mr} \le 1$, $r =1, \dots, R,\ m = 1, \dots, M-1$. As these constraints connect different levels of compliance table, the ILP problems can not be solved separately for each level, and we construct an ILP formulation for the whole table. 
Let $S_m$ be the event of having $m$ idle service engineers in the system. Then the objective function can be formulated as $\sum_{m=1}^M\mathbb{P}\left(S_m\right)\sum_{k=1}^K z_{mk}$. An accurate approximation of probabilities $\mathbb{P}\left(S_m\right)$ is provided in Appendix~\ref{apx:exp_cov}, equations~\eqref{eq:busy_m}. Finally, we are in position to provide the formulation for the MCRP compliance table:

\begin{minipage}{\linewidth-1cm}
\begin{flalign}
\label{mcrp_obj}
\max \quad &\sum_{m=1}^M\mathbb{P}\left(S_m\right)\sum_{k=1}^K z_{mk}\\
\label{mcrp_c1}
\text{s.t.} \quad & z_{mk} \le \sum_{r \in N_k}x_{mr}, &&  k = 1,\dots, \quad K, m = 1,\dots,M,\\
\label{mcrp_c2}
& \sum_{r=1}^{R}x_{mr} = m,  &&   m=1,\dots,M,\\
\label{mcrp_c3}
& x_{mr} - x_{m+1,r} \le \alpha_{mr}, &&  r = 1, \dots, R, \quad m = 1,\dots,M-1,\\
\label{mcrp_c4}
&\sum_{r=1}^R \alpha_{mr} \le 1,  &&  m = 1, \dots, M-1,\\
\label{mcrp_c5}
& \alpha_{mr} \ge 0, &&  r = 1, \dots, R, \quad m = 1,\dots,M-1,\\
\label{mcrp_c6}
& x_{mr} \in \{0, 1, \dots, m\} &&    r=1,\dots,R, \quad m=1,\dots,M,\\
\label{mcrp_c7}
& z_{mk} \in \{0,1\}  &&  k = 1, \dots K, \quad m=1,\dots,M.
\end{flalign}~
\end{minipage}


\subsection{MEXCRP compliance tables (RP3)} \label{sec:reloc:amexprep}

Consider a system where the number of idle service engineers is larger than the number required to cover all machines. In the MCRP approach presented in Section~\ref{sec:reloc:mcrp}, if each machine is covered by a service engineer, the location of the remaining service engineers does not affect the coverage. This may lead to inefficient allocation of service engineers. 
The hypothesis is that MEXCRP compliance tables introduced in~\cite{barneveld2016b} can solve this problem. In this algorithm, the main goal is to optimize expected covered demand, not the number of covered machines.

The problem can again be formulated as an ILP problem. 
Let the binary variable $y_{mki}$, $m = 1,\dots,M$, $k=1,\dots,K$, $i =1,\dots,m$, equal $1$ if and only if in the configuration for $m$ idle service engineers machine $k$ is covered by at least $i$ service engineers. Denote by $P_{mki}$, $m = 1,\dots,M$, $k=1,\dots,K$, $i =1,\dots,m$, the probability that in the configuration for $m$ idle service engineers a call from machine $k$ is responded to by the $i$th closest service engineer. For a given number of idle service engineers $m$, this probability can be approximated using equation~\eqref{eqn::ibusyprob} in Appendix~\ref{apx:exp_cov}. The MEXCRP compliance table formulation is as follows:

\begin{minipage}{\linewidth-1cm}
\begin{flalign}
\label{mexcrp_obj}
\max \quad &\sum_{m=1}^M\sum_{k=1}^K\sum_{i=1}^m P_{mki} y_{mki}\\
\label{mexcrp_c1}
\text{s.t.} \quad & \sum_{i = 1}^m y_{mki} \le \sum_{r \in N_k}x_{mr},  && k = 1,\dots, K, \quad  m = 1,\dots,M,\\
\label{mexcrp_c2}
& \sum_{r=1}^{R}x_{mr} \le m,  &&  m=1,\dots,M, \\
\label{mexcrp_c3}
& x_{mr} - x_{m+1,r} \le \alpha_{mr},  &&  r = 1, \dots, R, \quad  m = 1,\dots,M-1, \\
\label{mexcrp_c4}
&\sum_{r=1}^R \alpha_{mr} \le 1,   &&  m = 1, \dots, M-1, \\
\label{mexcrp_c5}
& \alpha_{mr} \ge 0,  &&   r = 1, \dots, R, \quad m = 1,\dots,M-1, \\
\label{mexcrp_c6}
& x_{mr} \in \{0, 1, \dots, m\} &&   r=1,\dots,R, \quad m=1,\dots,M, \\
\label{mexcrp_c7}
& y_{mki} \in \{0,1\} &&   k = 1, \dots, K, \quad m=1,\dots,M.
\end{flalign}~
\end{minipage}

\subsection{DMEXCLP heuristics (R4, R5)} \label{sec:reloc:dmexclp}

When compliance table relocation policy is used, the state that is achieved after relocation does not depend on the current state of the system, only on the number of idle service engineers. In contrast, DMEXCLP heuristic approach to relocation uses all the information about the current state to make a decision. This approach is more flexible than the compliance tables. We consider the DMEXCLP heuristic relocation policy introduced in~\cite{jagtenberg2015} and adjust it for our model.

There are two types of decision moments:
\begin{enumerate}
\item When a service is completed and there are no jobs assigned to the service engineer that just became idle. In this case, that service engineer must be dispatched to one of the base stations.
\item When an idle service engineer is dispatched to an incident, it should be decided whether one of the other idle service engineers should be relocated or not.
\end{enumerate}

According to the DMEXCLP relocation heuristic policy, the action that maximizes the expected covered demand is always chosen. In the first case, one service engineer is added sequentially to every base station, the expected covered demand is calculated, and the base station that leads to the best result is chosen. In the second case, all pairs of base stations $\left(r_1, r_2\right)$, with at least one service engineer at the base station $r_1$, are considered. We calculate the improvement in the expected covered demand after the relocation of a service engineer from station $r_1$ to station $r_2$. Suppose that $\left(r_1', r_2'\right)$ is the pair with the maximum improvement. If this improvement is positive, then we decide to relocate a service engineer from the station $r_1'$ to the station $r_2'$. If the maximum improvement is not positive, then no relocation happens. Note that for both situations the expected covered demand is computed only for the working machines.

The problem of large relocation times leading to the possibly poor performance can appear. In \cite{barneveld2018}, the authors impose restrictions on the relocations as a solution to this problem. There are three possible parameters that can be used to describe these restrictions. In the first type of decision moment the maximum relocation distance can be set. In this case, the best station is chosen among all stations within this distance from the machine where the service engineer is located. If there are no such base stations then the choice is made from all base stations. In the second type of decision moment the restriction can be imposed not only on the maximum relocation distance, but also on the minimum improvement in the expected covered demand. If the maximum relocation distance is set, than only the pairs $\left(r_1, r_2\right)$ with smaller distance are considered. If there are no such pairs, the relocation is forbidden. Note that this distance can differ from the distance for the first type of decision. If the minimum improvement threshold is set, then the relocation happens only if the improvement in expected covered demand exceeds this threshold. Setting this threshold equal to 0 means no restriction. The threshold larger than the number of machines leads to no relocation.

The optimal restriction parameters depend on the type of the system. However, there are no known results on how to find the optimal parameters for a given system. In the computational study in Section \ref{sec:num:reloc} below we consider this policy with different parameters and study the improvement that can be gained by parameter tuning.

\section{Numerical results}\label{sec:num}

In this section we present the setup and the results of our numerical experiments. To compare different policies we use simulation. In Section~\ref{sec:num:setup} we describe the type of the systems we use in our simulations, and the parameters defining the properties of the systems that affect the policies' performance. Section~\ref{sec:num:disp} presents the results comparing the dispatching policies introduced in Section~\ref{sec:disp} to each other with no relocation allowed. Next, in Section~\ref{sec:num:reloc} the relocation policies from Section~\ref{sec:reloc} are compared to each other, given a fixed dispatching policy. Finally, in Section~\ref{sec:num:opt} a heuristic policy based on the combination of the best dispatching and relocation policy from the previous sections is compared against the optimal policy for a small problem instance. We present the results of numerical experiments in tables for a limited set of parameter values. In the Appendix~\ref{apx:num} we provide extended tables for a wider range of systems.

\subsection{Setup of the numerical experiments}\label{sec:num:setup}

In this section we discuss the important parameters of the system, and how those parameters may affect the performance of various policies under consideration. Table \ref{tab::parameters} contains all parameters of the system. We now discuss the most important parameters and the relations between them in more detail.

\begin{table}[b]
\centering
\bgroup
\def\arraystretch{0.7}
\begin{tabular}{ll} \hline \hline
$K$ & number of machines \\
$R$ & number of base stations \\
$M$ & number of service engineers \\
& \\
$1/\lambda$ & average time till break down \\
$1/\mu$ & average repair time \\
$t^*$ & time limit \\
$d$ & map density \\
\hline \hline
\end{tabular}
\egroup
\caption{System parameters}
\label{tab::parameters}
\end{table}

\begin{figure}[t]
\centering
        \begin{subfigure}[b]{0.48\textwidth}
                \includegraphics[width=\linewidth]{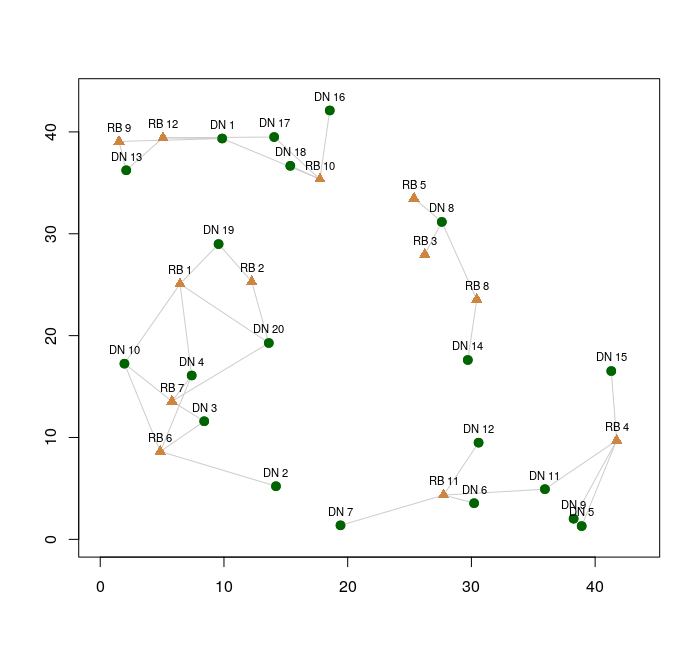}
                \caption{$d = 0.3$}
                \label{fig::density:0.3}
        \end{subfigure}%
        \begin{subfigure}[b]{0.48\textwidth}
                \includegraphics[width=\linewidth]{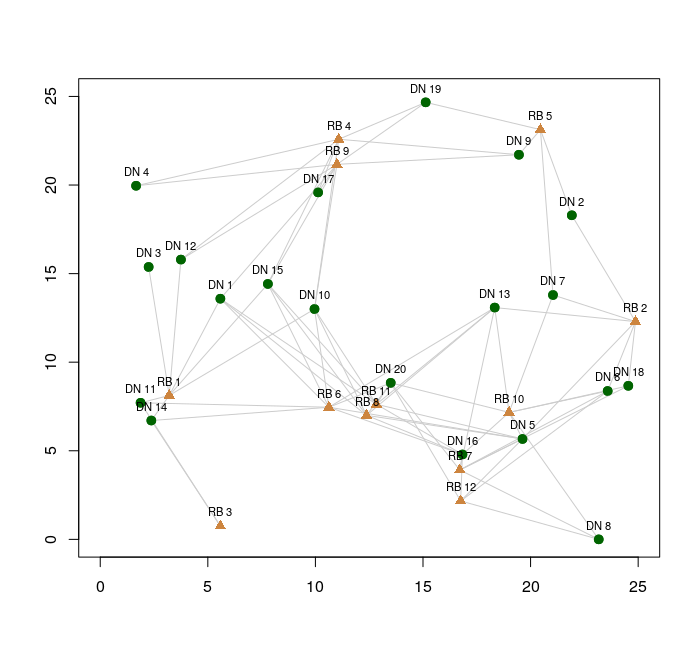}
                \caption{$d = 2$}
                \label{fig::density:2}
        \end{subfigure}%
        \caption{Maps of different density $d$, given $t^*=10$. Edges are shown for distances less than $t^*$.}\label{fig::density}
\end{figure}

\noindent
\textbf{Map density $\pmb{d}$.} Map density represents the relation between the average distance between nodes in the service region and the time limit $t^*$. The higher $d$ for the same value of $t^*$, the more dense the map, meaning that an average machine is covered by more base stations. Figure \ref{fig::density} provides an example with two randomly constructed maps for two different values of $d$, given the same value of $t^*$.

Map density can affect the optimal policy the following way. On one hand, if the map is sparse, the distances between base stations are large and relocation of service engineers may take undesirably long time. In the extreme case, each machine is covered by only one base station, so it is important that at least one service engineer is present at each station. On the other hand, in dense maps the distances between the machines are rather small compared to $t^*$, so when a machine breaks down it can be better to wait for a nearby busy service engineer to finish a repair of a machine instead of dispatching an idle service engineer.

\noindent
\textbf{Relation between $\pmb{\mu}$ and $\pmb{\lambda}$.} Given a fixed map and the number of service engineers, this relation influences the load of the system. If $\mu$ is increased with the fixed $\lambda$, the system becomes more loaded, that leads to decrease in the number of calls answered in time.

\subsection{Comparison of dispatching heuristics}\label{sec:num:disp}

In this section we show the performance of the policies described in Section~\ref{sec:disp}. To this end, we generate systems with different parameters and compare the fraction of calls answered in time under each of the five policies DP1-DP5. Relocation is not allowed, meaning that each service engineer is allocated to a fixed base station according to the model in Appendix~\ref{apx:opt_allocation}, and returns there after repairing a machine. The starting state for all policies is the same, chosen to maximize the expected covered demand.

We test the policies in a simulation. The number of machines, the number of bases and the failure rate are fixed at $K = 20$, $R = 12$ and $\lambda = 0.01$, respectively. We then change the number of engineers $M$, the map density $d$, the time limit $t^*$, and the the repair rate $\mu$. Those parameters in combination control the geographical structure of the service region and the load. For each combination of parameters we randomly generate 10 different maps, and run a simulation over a time horizon of 1000 time units. For our first experiment we compare the policies DP1, DP2 and DP3 used on the same maps. For each simulation run, we measure the fraction of calls responded to within the time limit $t^*$. Table \ref{tab:disp:nowaiting} contains the obtained results for a range of parameter values. One can see that the three policies perform very close to each other for all systems, and there is no policy that performs uniformly best. The performance decreases with load (the load increases for larger $M$ and lower $\mu$). We also observe that the performance of all three policies drops with the increase of $t^*$ for a given density $d$, especially dramatic for lower $d$ and lower load. Increasing the time limit $t^*$ for a fixed density $d$ results in the maps with larger distances compared to the average repair time.

\begin{table}[h]
\centering
\begin{tabular}{|l|l|l|l|l|l|l|l|l|}
\hline
\multicolumn{1}{|c|}{\multirow{2}{*}{\textbf{$\boldsymbol M$}}} & \multicolumn{1}{c|}{\multirow{2}{*}{\textbf{$\boldsymbol d$}}} & \multicolumn{1}{c|}{\multirow{2}{*}{\textbf{$\pmb{t^*}$}}} & \multicolumn{3}{c|}{\textbf{$\pmb{\mu = 0.2}$}}                                                              & \multicolumn{3}{c|}{\textbf{$\pmb{\mu = 0.05}$}}                                                             \\ \cline{4-9} 
\multicolumn{1}{|c|}{}                              & \multicolumn{1}{c|}{}                              & \multicolumn{1}{c|}{}                               & \multicolumn{1}{c|}{\textbf{DP1}} & \multicolumn{1}{c|}{\textbf{DP2}} & \multicolumn{1}{c|}{\textbf{DP3}} & \multicolumn{1}{c|}{\textbf{DP1}} & \multicolumn{1}{c|}{\textbf{DP2}} & \multicolumn{1}{c|}{\textbf{DP3}} \\ \hline
\multirow{4}{*}{10}                                 & \multirow{2}{*}{0.3}                               & 5                                                   & 0.92                             & 0.92                             & 0.93                             & 0.78                             & 0.79                             & 0.78                             \\ \cline{3-9} 
                                                    &                                                    & 20                                                  & 0.38                             & 0.36                             & 0.37                             & 0.29                             & 0.28                             & 0.28                             \\ \cline{2-9} 
                                                    & \multirow{2}{*}{2}                                 & 5                                                   & 0.99                             & 1.00                             & 1.00                             & 0.94 & 0.95                             & 0.95                             \\ \cline{3-9} 
                                                    &                                                    & 20                                                  & 0.87                             & 0.86                             & 0.86                             & 0.73 & 0.72                             & 0.73                             \\ \hline
\multirow{4}{*}{13}                                 & \multirow{2}{*}{0.3}                               & 5                                                   & 0.98                             & 0.98                             & 0.98                             & 0.92                             & 0.93                             & 0.94                             \\ \cline{3-9} 
                                                    &                                                    & 20                                                  & 0.84                             & 0.82                             & 0.82                             & 0.64                             & 0.60                             & 0.62                             \\ \cline{2-9} 
                                                    & \multirow{2}{*}{2}                                 & 5                                                   & 1.00                             & 1.00                             & 1.00                             & 0.99                             & 0.99                             & 0.99                             \\ \cline{3-9} 
                                                    &                                                    & 20                                                  & 0.98                             & 0.98                             & 0.98                             & 0.93                             & 0.94                             & 0.94                             \\ \hline
\end{tabular}
\caption{Fraction of calls answered in time for the policies DP1, DP2 and DP3.}
\label{tab:disp:nowaiting}
\end{table}

As the difference in performance between the first three policies is negligible, we then compare the simplest policy DP1 against the policies DP4 and DP5 in a separate experiment. Again, the three policies DP1, DP4 and DP5 are used in a simulation run on 10 randomly generated maps for each combination of parameter values. The obtained simulation results can be found in Table~\ref{tab:disp:waiting}. For most of the systems both policies with waiting outperform (or at least perform equal to) the traditional closest-first policy. The maximum relative improvement is over $50\%$. It increases with the decrease in map density and load, and increase in traveling times (remember that the traveling times increase with $t^*$ for a given $d$). In dense maps it is more likely that the machines are close to each other, and waiting can be beneficial only if the service times are relatively low compared to the traveling times. Comparing the performance of policies DP4 and DP5, one can see that including the remaining repair time leads to an improved performance, although marginal. Policy DP4, however, is probably more realistic, as it does not assume that the repair times are known up front.

\begin{table}[h]
\centering
\begin{tabular}{|l|l|l|l|l|l|l|l|l|}
\hline
\multicolumn{1}{|c|}{\multirow{2}{*}{\textbf{$\boldsymbol M$}}} & \multicolumn{1}{c|}{\multirow{2}{*}{\textbf{$\boldsymbol d$}}} & \multicolumn{1}{c|}{\multirow{2}{*}{\textbf{$\pmb{t^*}$}}} & \multicolumn{3}{c|}{\textbf{$\pmb{\mu = 0.2}$}}                                                              & \multicolumn{3}{c|}{\textbf{$\pmb{\mu = 0.05}$}}                                                             \\ \cline{4-9} 
\multicolumn{1}{|c|}{}                              & \multicolumn{1}{c|}{}                              & \multicolumn{1}{c|}{}                               & \multicolumn{1}{c|}{\textbf{DP1}} & \multicolumn{1}{c|}{\textbf{DP4}} & \multicolumn{1}{c|}{\textbf{DP5}} & \multicolumn{1}{c|}{\textbf{DP1}} & \multicolumn{1}{c|}{\textbf{DP4}} & \multicolumn{1}{c|}{\textbf{DP5}} \\ \hline
\multirow{4}{*}{10}                                 & \multirow{2}{*}{0.3}                               & 5                                                   & 0.92                             & 0.95                             & 0.96                             & 0.79                             & 0.81                             & 0.83                             \\ \cline{3-9} 
                                                    &                                                    & 20                                                  & 0.37                             & 0.83                             & 0.85                             & 0.31                             & 0.54                             & 0.64                             \\ \cline{2-9} 
                                                    & \multirow{2}{*}{2}                                 & 5                                                   & 0.99                             & 0.99                             & 0.99                             & 0.94                             & 0.95                             & 0.95                             \\ \cline{3-9} 
                                                    &                                                    & 20                                                  & 0.86                             & 0.97                             & 0.97                             & 0.72                             & 0.79                             & 0.88                             \\ \hline
\multirow{4}{*}{13}                                 & \multirow{2}{*}{0.3}                               & 5                                                   & 0.98                             & 0.98                             & 0.98                             & 0.92                             & 0.93                             & 0.92                             \\ \cline{3-9} 
                                                    &                                                    & 20                                                  & 0.77                             & 0.92                             & 0.94                             & 0.67                             & 0.79                             & 0.84                             \\ \cline{2-9} 
                                                    & \multirow{2}{*}{2}                                 & 5                                                   & 1.00                             & 0.99                             & 0.99                             & 0.99                             & 0.98                             & 0.99                             \\ \cline{3-9} 
                                                    &                                                    & 20                                                  & 0.98                             & 0.95                             & 0.98                             & 0.94                             & 0.93                             & 0.96                             \\ \hline
\end{tabular}
\caption{Fraction of calls answered in time for the policies DP1, DP4 and DP5.}
\label{tab:disp:waiting}
\end{table}

\subsection{Comparison of relocation heuristics}\label{sec:num:reloc}

In this section we present simulation results comparing different relocation policies. The dispatching policy is fixed to the policy DP4 as one of the best performing policies from Section~\ref{sec:num:disp}.

As mentioned in Section~\ref{sec:reloc:dmexclp} above, there are three parameters that define restrictions for the policy RP5. Those are the maximum relocation distance upon redeployment after repair is finished, the maximum relocation distance upon dispatching, and the minimum performance improvement for relocation upon dispatching. To fit these parameters for the RP5 policy, we simulated the system for the first and the second parameters equal to $0.5t^*$, $t^*$, $2t^*$ and $100t^*$, and the third parameter equal to 0, 1, 5 and 100. The best result for each type of the system was chosen and used as an input for the policy RP5.

The relocation policies are compared using simulation. We set $K = 20$, $R = 10$ and $\lambda = 0.01$. This time we also fix the number of service engineers equal to $M=13$. For each type of the system we generate $10$ random maps, run simulation for each of the maps, and measure the fraction of calls responded to within the time limit for each of the five policies. The results can be found in Table \ref{tab:reloc}. One can see that compliance tables (RP2 and RP3) demonstrate poor performance for most systems. The only type of the systems for which RP2 and RP3 perform better than the static policy RP1 is those with high density and repair times much larger than $t^*$. The reason is that both MCRP and MEXCRP algorithms ignore the distances, so when the distances are large, those policies lead to inefficient relocations and poor performance.

\begin{table}[t]
\centering
\begin{tabular}{|l|l|l|l|l|l|l|l|l|l|l|l|}
\hline
\multicolumn{1}{|c|}{\multirow{2}{*}{\pmb{$d$}}} & \multicolumn{1}{c|}{\multirow{2}{*}{\pmb{$t^*$}}} & \multicolumn{5}{c|}{\pmb{$\mu=0.2$}}                                                                                                                                           & \multicolumn{5}{c|}{\pmb{$\mu=0.05$}}                                                                                                                                          \\ \cline{3-12} 
\multicolumn{1}{|c|}{}                              & \multicolumn{1}{c|}{}                               & \multicolumn{1}{c|}{\textbf{RP1}} & \multicolumn{1}{c|}{\textbf{RP2}} & \multicolumn{1}{c|}{\textbf{RP3}} & \multicolumn{1}{c|}{\textbf{RP4}} & \multicolumn{1}{c|}{\textbf{RP5}} & \multicolumn{1}{c|}{\textbf{RP1}} & \multicolumn{1}{c|}{\textbf{RP2}} & \multicolumn{1}{c|}{\textbf{RP3}} & \multicolumn{1}{c|}{\textbf{RP4}} & \multicolumn{1}{c|}{\textbf{RP5}} \\ \hline
\multirow{2}{*}{0.3}                                & 5                                                   & 0.93                              & 0.69                              & 0.71                              & 0.89                              & 0.98                              & 0.76                              & 0.70                              & 0.70                              & 0.85                              & 0.95                              \\ \cline{2-12} 
                                                    & 20                                                  & 0.77                              & 0.43                              & 0.43                              & 0.48                              & 0.92                              & 0.52                              & 0.25                              & 0.25                              & 0.43                              & 0.82                              \\ \hline
\multirow{2}{*}{2}                                  & 5                                                   & 0.97                              & 0.87                              & 0.89                              & 0.99                              & 0.99                              & 0.87                              & 0.92                              & 0.96                              & 0.98                              & 0.99                              \\ \cline{2-12} 
                                                    & 20                                                  & 0.92                              & 0.74                              & 0.79                              & 0.79                              & 0.97                              & 0.64                              & 0.49                              & 0.55                              & 0.75                              & 0.96                              \\ \hline
\end{tabular}
\caption{Fraction of calls answered in time for the policies RP1-RP5 combined with DP4.}
\label{tab:reloc}
\end{table}

Policy RP4 outperforms compliance tables for all systems because it uses more information about the state of the system to make relocation decisions. However, it also ignores the distances, so for maps with large distances and small density we observe that it performs worse than the static policy. Finally, policy RP5 leads to good results for all of the systems. The fraction of calls answered in time stays above $80\%$, even for the systems with high load, where all other policies result in less than $60\%$ of calls answered in time. The difference in performance between RP5 and RP4 shows the importance of relocation restrictions and accurate tuning of the parameters of these restrictions.

\subsection{Optimal policy performance}\label{sec:num:opt}

\begin{figure}
\centering
\includegraphics[width = 0.6\textwidth]{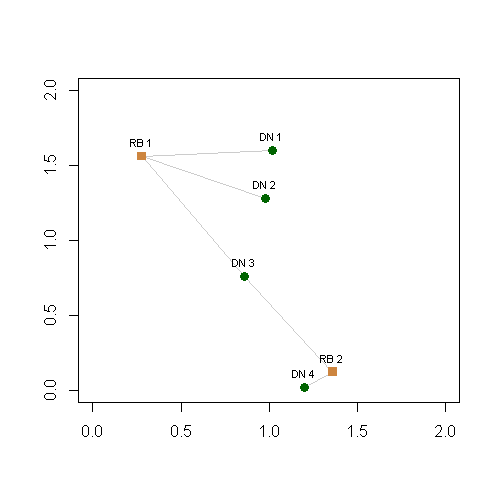}
\captionsetup{width = \textwidth}
\caption{Sample map with two base stations (RB) and four machines (DN).}
\label{fig::mapmdp}
\end{figure}

In this section, we benchmark the best performing heuristic against the optimal policy. Heuristic policy is combined of the dispatching policy DP4 and the relocation policy RP5 (i.e., minimal response time + DMEXCLP). To obtain the optimal policy, we use the discrete-time model described in Appendix~\ref{apx:dt_model}. Due to computational complexity of finding the optimal policy, we do this only for one small system depicted on Figure \ref{fig::mapmdp}. There are two base stations, two service engineers, and four machines. The failure rate is $\lambda = 0.01$ and the repair rate is $1$. The time limit is set to $t^* = 3$.

We use policy iteration to derive the optimal policy. Then we run $10$ iterations of simulation under the optimal policy and under the heuristic. The obtained average fraction of calls answered in time is $0.87$ for the optimal policy and $0.82$ under the DP4+RP5 heuristic policy, which is $5.7\%$ less than the optimal performance. We see that the heuristic performs close to the optimal at least for this instance. Unlike the optimal policy, however, it can be easily derived for real-life systems with significantly larger state spaces.

\section{Conclusion}\label{sec:conclusion}

In this paper, we studied the problem of real-time management of service engineers. We drew inspiration from the vast research in EMS domain to develop a number of scalable heuristics that lead to a low number of late arrivals to the emergency calls. We compared the performance of multiple heuristics, examining if there is any that works well irrespective of the network structure. We conducted extensive computational experiments where a range of dispatching and relocation policies were compared against each other for various types of systems. One of the best combined policies was then tested in a simulation against the optimal policy showing close to the optimal performance.

The dispatching policies were compared with the relocation policy fixed to a static one, where idle service engineers reside at preassigned base stations. Five dispatching policies were compared in a simulation for systems with different parameters. Parameter defining the system include the number of service engineers, the repair rate, the time limit and the map density. The policy that assigns the service engineers with the smallest response time outperformed other policies for most of the systems we considered. We showed that it may be beneficial not to dispatch an idle service engineer immediately upon a break down of a machine, but wait for another service engineer to finish repair. We also demonstrated that accurate estimation of repair time can lead to improvement in performance, although only for the systems with high load.

The performance of five relocation policies was measured with the fixed dispatching policy that always chooses the service engineer with the minimal response time. For most of the systems, compliance tables performed worse than the closest-first policy, except for those where the distances are small compared to the average repair time. The numerical results favor the DMEXCLP relocation policy. However, its performance depends on the choice of the restriction parameters, such as the maximum distance of relocation and the minimum improvement in the expected covered demand. Without these restrictions it performs worse than the closest-first policy for systems with large distances. However, when the restriction parameters are carefully chosen, it outperforms other policies by up to $60\%$ in fraction of calls answered in time.

\noindent
\textbf{Acknowledgements}

This research was partly funded by an NWO grant, under contract number 438-15-506.

\bibliographystyle{plain}
\bibliography{references}

\newpage
\appendix

\section{Discrete-time Model}\label{apx:dt_model}

In this section, we construct a discrete-time approximation of the original process formulated in Section~\ref{sec:model}. To that end, we discretize time and the service region. The service region is approximated by the one where all distances are rounded to the integer numbers. Recall that for the continuous time model, we look at the state of the process at each moment an event happens. In the discrete-time model we look at the state of the process at each time point, where several events can happen between subsequent time points. The resulting process has a finite state space that allows us to perform the policy iteration algorithm. In practice, the length of the time unit is an important decision to make. Small duration of the time unit provides good approximation of the real process, however it may also lead to a large state space that is computationally intractable.

We then define the state space based on the continues-time model, making sure it is finite. For the continuous-time process the state is described by a tuple $s = \left(t, e, \mathfrak{m}, \kappa\right)$, where $t$ is the time, $e$ is the event, $\mathfrak{m}$ is the state of the service engineers and $\kappa$ is the state of the machines. For the discrete-time process, we look at the state of the process at each time unit. So the time in state $s_n$ is deterministic $t(s_n) = t_n = n$, and the optimal action in state $s$ should not depend on time $t(s)$. We then omit the time component of state in the discrete-time version of the process. The vector $\mathfrak{m}$ contains the pairs of destinations and distances to those destinations for each service engineer. For the discrete-time model, however, distances take only integer values, so there is only a finite number of possible locations of service engineers. Vector $\kappa$ contains the state of each of the machines: the waiting time if the machine is broken, $0$ if it is working, and $-1$ if it is in repair. To have a finite number of possible states for each machine, the waiting time is bounded by the time limit $t^*$. If the waiting time of a broken machine $k$ reaches $t^*$, its state remains unchanged ($\kappa_k = t^*$) until repair starts.

In the discrete-time process more than one event can happen during the transition from state $s_n$ to state $s_{n+1}$ (for example, two machines may break down). Hence, we define the event $e$ as a triplet of sets $e = \left(\mathcal{K}_1, \mathcal{K}_2, \mathcal{M}_a\right)$, where $\mathcal{K}_1$ is the set of machines that got broken during the last time unit, $\mathcal{K}_2$ is the set of machines that got repaired during the last time unit, and $\mathcal{M}_a$ is the set of service engineers that arrived to their destinations during the last time unit. It is possible that no events happened during the transition, and all three sets are empty. 

\subsection{State transitions} \label{apx:dt_model:trans}

To determine transition probabilities, we basically need to derive the probability of a certain event $e = \left(\mathcal{K}_1, \mathcal{K}_2, \mathcal{M}_a\right)$ happening in a given state $s$. This probability depends only on the first two components, as the last component $\mathcal{M}_a$ is deterministic and depends only on the previous state of the system. The probability that a working machine will break down during one time unit is $p = 1-e^{-\lambda}$, and the probability that a repair in progress will end during one time unit is $q = 1 - e^{-\mu}$. Denote by $W(s)$ the number of all working machines in state $s$, and by $H(s)$ the number of machines in repair. Then the probability of event $e = \left(\mathcal{K}_1, \mathcal{K}_2, \mathcal{M}_a\right)$ happening in state $s$ is equal to
$$P(s, e)=P(s, \mathcal{K}_1, \mathcal{K}_2) = p^{|\mathcal{K}_1|}(1-p)^{W(s)-|\mathcal{K}_1|} q^{|\mathcal{K}_2|}(1-q)^{H(s)-|\mathcal{K}_2|}.$$

The next state $s_{n+1}$ of the process depends only on the current state $s_n$ of the process, the action $a_n$ taken, and the random components of the event ($\mathcal{K}_1$ and $\mathcal{K}_2$):
$$s_{n+1} = \Phi\left(s_n, a_n, \mathcal{K}_1, \mathcal{K}_2\right).$$

\subsection{Action space} \label{apx:dt_model:actions}

Let $\mathcal{F}\left(s\right) = \{m \in 1,\dots,M \mid l_m(s) \in \mathcal{R}\}$ denote the set of all idle service engineers, and $Q(s)$ the set of all machines in the queue in state $s$. Note that the set $\mathcal{F}\left(s\right)$ includes the set $\mathcal{M}_a$, and the set $Q(s)$ includes the set $\mathcal{K}_1$. An action $a$ in state $s$ consists of three binary matrices $\pmb{X}$, $\pmb{Y}$ and $\pmb{Z}$. Matrix $\pmb{X}$ describes the dispatching decision for all machines in the set $Q(s)$. Let $X_{mk}=1$ if service engineer $m$ is dispatched to machine $k$, and $X_{mk}=0$ otherwise. Matrix $\pmb{Y}$ describes the redeployment decision for the service engineers that finished repairing machines in the set $\mathcal{K}_2$. Let $Y_{ml}=1$ if service engineers $m$ is redeployed to location $l$, and $Y_{ml}=0$ otherwise. Matrix $\pmb{Z}$ describes the relocation decision. Let $Z_{mr} = 1$ if service engineer $m$ is relocated to base station $r$, and $Z_{mr} = 0$ otherwise. The action space in state $s$ is given by
\begin{align}
\nonumber \mathcal{A}(s) = \bigg\{(\pmb{X}, \pmb{Y}, \pmb{Z}) \mid &\sum_{m \in \mathcal{F}(s)}X_{mk} + \sum_{m \in \mathcal{M}:l_m(s)\in \mathcal{K}_2}Y_{mk} \le 1 \ \forall k \in Q(s); \\
\nonumber &\sum_{l \in Q(s)\cup\mathcal{R}}Y_{ml} = 1 \ \forall m \in \mathcal{M}:l_m(s)\in \mathcal{K}_2;  \\ 
\nonumber &\sum_{\substack{m \in \mathcal{F}(s),\\ r \in \mathcal{R}}}Z_{mr}\le 1; \ \sum_{\substack{m \in \mathcal{F}(s),\\ r \in \mathcal{R},\\ k \in Q(s)}} X_{mk}Z_{mr}=0; \\
& \mathbb{I}(\sum{\substack{m \in \mathcal{F}(s),\\ k \in \mathcal{K}_1}} X_{mk}=0) \mathbb{I}(\mathcal{K}_2=\emptyset) \sum_{\substack{m \in \mathcal{F}(s),\\ r \in \mathcal{R}}}Z_{mr} = 0 \bigg\}, \label{eqn:A}
\end{align}
where the constraints ensure that not more than one service engineer is dispatched to each broken machine, each service engineer that finished a job is redeployed to exactly one location, at most one relocation is made and the relocated service engineer differs from the dispatched ones, relocation is done only if at least one service engineer was dispatched to one of the newly broken machines or at least one service engineer finished a repair.


\subsection{Costs} \label{apx:dt_model:costs}
We apply the following cost structure. If a call arrives from machine $k$ and a service engineer does not reach this machine within the time limit $t^*$, a cost $1$ is incurred. Moreover, a small cost $0<\epsilon \ll 1$ is incurred per time unit of waiting for service over $t^*$. This second penalty is used to ensure dispatching is done by the optimal policy. Our goal is to find actions that minimize the long-run discounted penalty. All travel costs and other operational costs are ignored, but could be readily added to the model.

Denote by $c(s_n, a_n, s_{n+1})$ the costs incurred during the transition period from state $s_n$ to state $s_{n+1}$ when action $a_n$ is taken. The first component is a unit penalty for each machine who's waiting time exceeds $t^*$ during transition. The extra costs are incurred for the total waiting time over $t^*$. Then in total
$$c(s_n, a_n, s_{n+1}) = |\left\{k \in \mathcal{K} \mid \kappa_k(s_n) < t^* \text{ and } \kappa_k(s_{n+1}) = t^*\right\}| + \epsilon |\left\{k \in \mathcal{K} \mid \kappa_k(s_n) = t^*\right\}|.$$
It is important to note that our goal is to maximize the fraction of calls answered in time. The penalty $\epsilon$ is introduced only to prevent the situation of leaving some machines broken forever, which would otherwise be optimal. So $\epsilon$ is set to be small. In the computational experiments, we set $\epsilon = 0.001$, so it does not affect the optimal policy.

\subsection{Optimality equations} \label{apx:dt_model:opt_eq}

Consider a discounted version of the process $\{s_n, \ n=0,1,\dots\}$ with discount factor $\gamma$~\cite{puterman2014}. Denote by $V_{\pi}(s)$ the expected total discounted costs under policy $\pi$ when starting in state $s$:

\begin{equation*}
V_{\pi}(s) = \mathbb{E}\left[\sum_{n=0}^\infty \gamma^{n}c(s_n, \pi(s_n), s_{n+1}) \mid s_0=s \right].
\end{equation*}
If policy $\pi^*$ is the optimal policy, then $V_{\pi^*}(s)$ satisfies the Bellman optimality equation:
\begin{equation*}
V_{\pi^*}(s) = \min_{a \in \mathcal{A}(s)}\left\{\mathbb{E}_a\left[c(s, a, s') + \gamma V_{\pi^*}(s')\right]\right\}, \quad \forall s \in S,
\end{equation*}
where $s' = \Phi\left(s, \pi(s), \omega(s, a)\right)$ is the next state of the process when action $a$ is taken in state $s$, and
\begin{equation*}
\pi^*(s) = \argmin_{a \in \mathcal{A}(s)}\left\{\mathbb{E}_a\left[c(s, a, s') + \gamma V_{\pi}(s')\right]\right\}, \quad \forall s \in S.
\end{equation*}
\noindent
Due to the curse of dimensionality, finding the optimal policy $\pi^*$ is computationally intractable for realistic-sized systems. Hence in Sections \ref{sec:disp} and \ref{sec:reloc} we focus on various scalable heuristic approaches to the above problem.

The state of the discrete-time process is represented by a triplet $s = \left(e, \mathfrak{m}, \kappa\right)$. As the set of possible states of service engineers and machines, and the set of possible events are all finite, the state space of the process is finite. As the transitions and the costs depend only on the current state and the random component, the resulting process is a finite-state Markov decision process.

\section{Expected covered demand approximation}\label{apx:exp_cov}

In this section we derive an approximation of the expected covered demand, where we follow the procedure introduced by Larson \cite{larson1974, larson1975}, and apply it to the model described in Section~\ref{sec:model}. Given the locations of the service engineers, the expected covered demand estimates the long-term fraction of calls that will be answered in time. We use this metric for several dispatching and relocation policies (see Sections \ref{sec:disp} and \ref{sec:reloc}).

First, we consider the process $C = \left\{C_n, \ n = 1, 2, \dots \right\}$, that approximates the number of broken machines in the $n^{th}$ state of the original process $s_n$, $n=1,2,\dots$, where $C_n = |\left\{k\in\mathcal{K} \mid \kappa_k(s_n) \ne 0 \text{ or } e = \text{"machine $k$ breaks down"}\right\}|$. We compute the steady-state distribution of this process. 

The time till breakdown of a machine is exponentially distributed with rate $\lambda$.  The time a machine stays broken is also exponentially distributed with rate $\hat{\mu}$. This time includes the traveling time of a service engineers and the duration of the repair. If all service engineers are busy, and the machine is put in the queue first, then the waiting time in the queue is not included.

\begin{figure}[b]
\centering
\includegraphics[width = \textwidth]{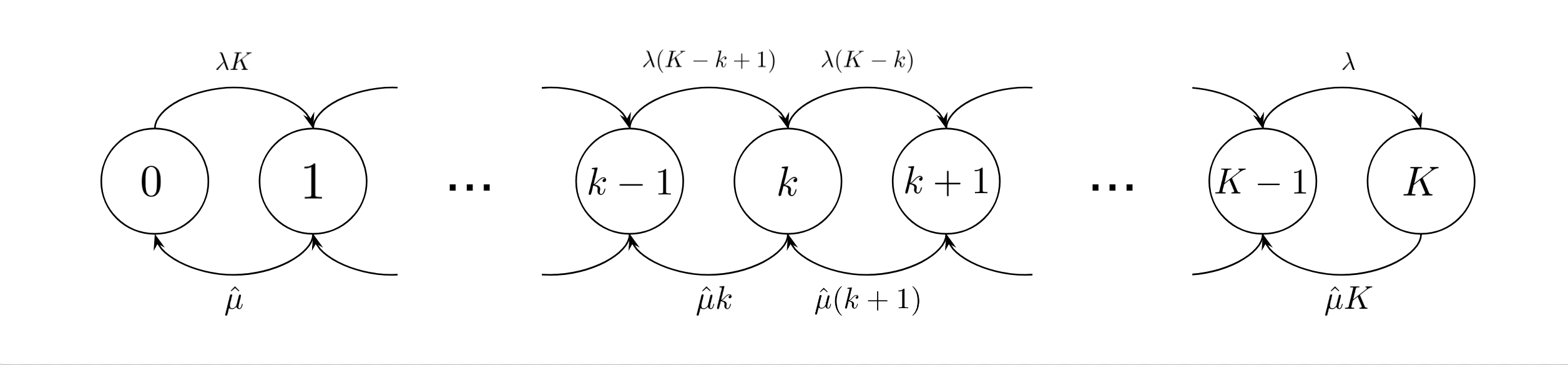}
\caption{State diagram of the discrete-time process $C$.}
\label{fig:process_c}
\end{figure}

The process $C$ is a Markov process. The state space of the process is $\{0,1,\dots,K\}$ (see Figure~\ref{fig:process_c}), so it is a finite-state process. There are two possible transitions from a state $C_n$ with $k$ broken machines:
\begin{itemize}
\item The event in state $s_{n+1}$ is of type {\it "a machine breaks down"}. Then $C_{n+1} = k+1$. The rate of this transition is equal to $\lambda(K-k)$. (This transition is not possible if $k=K$.)
\item The event in state $s_{n+1}$ is of type {\it "a repair ends"}. Then $C_{n+1} = k-1$. The rate of this transition $\hat{\mu}\cdot\text{\# machines in repair} = \hat{\mu}\cdot \min\{k, M\}$. (This transition is not possible if $k=0$.)
\end{itemize}

Let $P(k)$ denote the stationary probability of being in state $k$. Then the balance equations for $C_n$ can be formulated as follows:
\begin{equation*}
\begin{cases}
 \quad\lambda K P(0) = \hat{\mu} P(1), \\
 \quad\left(\lambda (K-k)+\hat{\mu} k\right) P(k) = \lambda(K-k+1) P(k-1) + \hat{\mu}(k+1) P(k+1), \quad & k = 2, \dots, M-1, \\
\quad \left(\lambda (K-k)+\hat{\mu} M\right) P(k) = \lambda(K-k+1) P(k-1) + \hat{\mu} M P(k+1), \qquad & k = M, \dots, K-1, \\
\quad \hat{\mu}M P(K) = \lambda P(K-1).
\end{cases}
\end{equation*}
\noindent
One can check that 
\begin{equation}
P(k)= \begin{cases} 
\binom{K}{k} \left(\frac{\lambda}{\hat{\mu}}\right)^k P(0), \qquad &k = 0,\dots, M-1, \\ 
\frac{k!}{M!M^{k-M}} \binom{K}{k} \left(\frac{\lambda}{\hat{\mu}}\right)^k P(0), \qquad &k = M,\dots, K
\end{cases}
\end{equation}
is the solution of the balance equations. Adding the normalization equation $\sum_{k=0}^K P(k) = 1$ to the system we get
\begin{equation}
P(0) = \left[\sum_{k=0}^{M-1}\binom{K}{k} \left(\frac{\lambda}{\hat{\mu}}\right)^k + \sum_{k=M}^{K} \frac{k!}{M!M^{k-M}} \binom{K}{k} \left(\frac{\lambda}{\hat{\mu}}\right)^k \right]^{-1} .
\end{equation}
\noindent
Using these formulas one can calculate $P(k)$ for $k=0,1, \dots, K$.

Now, let $S_m$ be the event of having $m$ busy service engineers, then
\begin{equation}\label{eq:busy_m}
\mathbb{P}(S_m) = \begin{cases}
P(m), \qquad &m = 1, \dots, M-1, \\
\sum_{k=M}^K P(k), \qquad &m = M.
\end{cases}
\end{equation}

We consider a system where relocation is not allowed. It means that each service engineer is assigned to a base station, and returns to that base station after each repair. Suppose that the system is in state $s$ where all service engineers are at their base stations. Assume that the dispatching policy is fixed per machine, such that if machine $k$ breaks down, we first send the service engineer $m_{1}^{(k)}$. Then, if that service engineer is busy, $m_{2}^{(k)}$ is dispatched, and so on. Assume also that first we try to send the service engineers that can reach the machine within the time limit, and, if all of them are busy, the service engineers that cannot arrive on time. An example of such dispatching policy is the closest-first policy that always dispatches the closest service engineer.

Let us compute the probability that service engineer $m_{i}^{(k)}$ is dispatched to machine $k$, meaning that all service engineers $m_{1}^{(k)},\dots,m_{i-1}^{(k)}$ are busy. If $B_i$ is the event that service engineer $m_{i}^{(k)}$ is busy, $F_i$ is the event that service engineer $m_{i}^{(k)}$ is idle, and $S_m$ is the event that there are $m$ busy service engineers in the system, then
\begin{multline}
\mathbb{P}\left(B_1\dots B_{i-1}F_i\right) = \sum_{m=i}^M \mathbb{P}\left(B_1\dots B_{i-1}F_i\right|S_m)\mathbb{P}(S_m) \\= \sum_{m=i}^M \mathbb{P}(S_m)\mathbb{P}\left(F_i\right|S_mB_1\dots B_{i-1})\mathbb{P}\left(B_{i-1}\right|S_mB_1\dots B_{i-2})\dots\mathbb{P}\left(B_1|S_m\right).
\end{multline}

The probabilities $P(S_m)$ can be computed using equations~\eqref{eq:busy_m}. Other terms can be approximated by assuming that all service engineers have the same load and are independent of each other. Under these assumptions we get
\begin{equation} 
\begin{gathered}
\mathbb{P}\left(B_1|S_m\right) = \frac{m}{M} \\
\mathbb{P}\left(B_2|S_m B_1\right) = \frac{m-1}{M-1} \\
\vdots \\
\mathbb{P}\left(F_i\right|S_mB_1\dots B_{i-1}) = 1 - \frac{m-i+1}{M-i+1} = \frac{M-m}{M-i+1}.
\end{gathered}
\end{equation}

Finally, we can approximate 
\begin{multline}
\mathbb{P}\left(B_1\dots B_{i-1}F_i\right) = \sum_{m=i-1}^M \mathbb{P}(S_m)\mathbb{P}\left(F_i\right|S_mB_1\dots B_{i-1})\mathbb{P}\left(B_{i-1}\right|S_mB_1\dots B_{i-2})\dots\mathbb{P}\left(B_1|S_m\right) \\ \approx \sum_{m=i-1}^M \mathbb{P}(S_m)\cdot \frac{M-m}{M-i+1} \cdot \dots \cdot \frac{m}{M} = \sum_{m=i-1}^M \left(M-m\right)\mathbb{P}(S_m)\cdot\frac{m!(M-i)!}{(m-i+1)!M!}.
\end{multline}

Let us denote 
\begin{equation} \label{eqn::ibusyprob}
P_i = \sum_{m=i-1}^M \left(M-m\right)\mathbb{P}(S_m)\cdot\frac{m!(M-i)!}{(m-i+1)!M!}, \quad i=1,\dots,M.
\end{equation}

Let the binary variable $z_{ki}$ equal $1$ if the $i^\text{th}$ closest service engineer to machine $k$ can reach it in time. The probability that a call from machine $k$ will be answered in time is $\sum_{i=1}^{M}P_i z_{ki}$, and the expected covered demand can be approximated by
\begin{equation}
\frac{1}{K}\sum_{k=1}^K\sum_{i=1}^{M}P_i z_{ki}.
\label{eqn::expcovdem}
\end{equation}

Note that for a given system all parameters needed to calculate expression \eqref{eqn::expcovdem} are known, except for the parameter $\hat{\mu}$. This parameter is hard to calculate in practice as it depends on the policy. For computational study, we first assume $\hat{\mu} = 1/(t^*+1/\mu)$ and then run several iterations of simulation to find a better approximation for $\hat{\mu}$.

\section{Optimal allocation of service engineers}\label{apx:opt_allocation}

In this section we formulate an integer linear programming problem that finds the optimal allocation of the service engineers to the base stations maximizing the expected covered demand.

Let the decision variables $x_r$, $r = 1, \dots, R$, represent the number of service engineers at base station $r$, and $z_{ki}$ indicate if service engineer $m_{i}^{(k)}$ can reach machine $k$ in time. The objective function is the expected covered demand that we approximate with \eqref{eqn::expcovdem}.

The total number of service engineers is $M$, so $\sum_{r=1}^{R}x_r = M$. The variables $z_{ki}$, $k = 1,\dots, K$, $i = 1,\dots,M$, and the variables $x_r$, $r=1,\dots,R$, are connected by the equation $\sum_{i=1}^{M}z_{ki} = \sum_{r \in N_k}x_r, \ k = 1,\dots, K$, where $N_k$ is the set of all bases from which machine $k$ can be reached in time. The problem can be formulated as follows:

\begin{equation}
\begin{aligned}
\max \quad &\sum_{k=1}^K\sum_{i=1}^{M}P_i z_{ki} \\
\text{s.t.} \quad &\sum_{i=1}^{M}z_{ki} \le \sum_{r \in N_k}x_r, && k = 1,\dots, K, \\
& \sum_{r=1}^{R}x_r \le M, \\
& x_r = 0, 1, 2, \dots &&  r = 1, \dots R, \\
& z_{ki} \in \{0, 1\} && k = 1,\dots,K, \quad r = 1, \dots, R,
\end{aligned}
\label{eqn::expcovdemalloc}
\end{equation}
where constraints are relaxed with inequalities. The total number of decision variables is $R+KM$ and the total number of constraints equals $K+1$.

\newpage
\section{Extended Computational Results}\label{apx:num}

\begin{table}[ht]
\small
\centering
\bgroup
\def\arraystretch{0.65}
\begin{tabu}to \textwidth {|X[c]|X[c]|X[c]|X[c]X[c]X[c]|X[c]X[c]X[c]|X[c]X[c]X[c]|X[c]X[c]X[c]|} 
\hline
\multirow{2}{*}{$M$}& \multirow{2}{*}{$d$} &\multirow{2}{*}{$t^*$} &\multicolumn{3}{c|}{$\mu = 0.2$}&\multicolumn{3}{c|}{$\mu = 0.1$}&\multicolumn{3}{c|}{$\mu = 0.05$}&\multicolumn{3}{c|}{$\mu = 0.02$}\\ \cline{4-15}
&&& DP1 &DP2 & DP3 & DP1 &DP2 & DP3 & DP1 &DP2 & DP3 & DP1 &DP2 & DP3 \\ \hline
\multirow{12}{*}{10}&\multirow{4}{*}{0.3}&5&0.92&0.92&0.93&0.88&0.88&0.89&0.78&0.79&0.78&0.44&0.44&0.42\\
&&10&0.81&0.79&0.79&0.76&0.73&0.74&0.60&0.57&0.57&0.33&0.33&0.31\\
&&20&0.38&0.36&0.37&0.32&0.32&0.32&0.29&0.28&0.28&0.25&0.25&0.24\\
&& 50&0.23&0.24&0.24&0.26&0.26&0.26&0.22&0.23&0.23&0.22&0.22&0.21\\ \cline{2-15}
&\multirow{4}{*}{1}& 5&0.96&0.97&0.97&0.94&0.95&0.95&0.89&0.90&0.89&0.61&0.62&0.60\\
&&10&0.91&0.92&0.91&0.86&0.87&0.87&0.76&0.76&0.76&0.44&0.45&0.42\\
&&20&0.59&0.56&0.57&0.54&0.54&0.53&0.45&0.45&0.45&0.29&0.30&0.29\\
&&50&0.34&0.34&0.35&0.34&0.34&0.34&0.37&0.37&0.37&0.28&0.29&0.28\\ \cline{2-15}
&\multirow{4}{*}{2}& 5&0.99&1.00&1.00&0.98&0.98&0.98&0.94&0.95&0.95&0.75&0.75&0.73\\
&&10&0.97&0.97&0.98&0.95&0.95&0.96&0.90&0.90&0.90&0.62&0.63&0.61\\
&& 20&0.87&0.86&0.86&0.81&0.82&0.81&0.73&0.72&0.73&0.47&0.47&0.46\\
&&50&0.53&0.52&0.52&0.49&0.49&0.50&0.50&0.50&0.50&0.45&0.45&0.44\\ \hline
\multirow{12}{*}{13}&\multirow{4}{*}{0.3}& 5&0.98&0.98&0.98&0.96&0.97&0.97&0.92&0.93&0.94&0.80&0.81&0.79\\
&&10&0.94&0.94&0.95&0.92&0.93&0.93&0.88&0.89&0.89&0.68&0.68&0.67\\
&&20&0.84&0.82&0.82&0.73&0.68&0.69&0.64&0.60&0.62&0.44&0.43&0.42\\
&& 50&0.32&0.33&0.33&0.32&0.32&0.32&0.29&0.30&0.30&0.28&0.29&0.28\\ \cline{2-15}
&\multirow{4}{*}{1}&5&0.99&1.00&1.00&0.98&0.99&0.99&0.97&0.97&0.97&0.86&0.88&0.87\\
&&10&0.98&0.98&0.99&0.96&0.96&0.97&0.94&0.95&0.95&0.80&0.82&0.80\\
&&20&0.90&0.89&0.91&0.91&0.91&0.91&0.84&0.83&0.83&0.62&0.61&0.60\\
&&50&0.44&0.45&0.45&0.43&0.43&0.43&0.44&0.45&0.45&0.40&0.40&0.40\\ \cline{2-15}
&\multirow{4}{*}{2}& 5&1.00&1.00&1.00&1.00&1.00&1.00&0.99&0.99&0.99&0.93&0.94&0.94\\
&& 10&1.00&1.00&1.00&0.99&0.99&0.99&0.98&0.98&0.98&0.90&0.91&0.91\\
&&20&0.98&0.98&0.98&0.96&0.96&0.97&0.93&0.94&0.94&0.81&0.79&0.79\\
&&50&0.65&0.63&0.63&0.60&0.60&0.60&0.58&0.56&0.56&0.58&0.58&0.58\\ \hline
\multirow{12}{*}{16}&\multirow{4}{*}{0.3}& 5&0.99&0.99&0.99&0.99&0.99&0.99&0.97&0.98&0.98&0.93&0.94&0.94\\
&&10&0.98&0.99&0.99&0.96&0.97&0.97&0.95&0.96&0.96&0.90&0.91&0.91\\
&& 20&0.94&0.95&0.96&0.92&0.92&0.93&0.91&0.90&0.91&0.81&0.81&0.80\\
&& 50&0.51&0.51&0.51&0.52&0.50&0.53&0.48&0.46&0.48&0.41&0.41&0.41\\ \cline{2-15}
&\multirow{4}{*}{1}& 5&1.00&1.00&1.00&0.99&1.00&1.00&0.99&0.99&1.00&0.96&0.97&0.97\\
&&10&0.99&0.99&0.99&0.99&0.99&0.99&0.98&0.98&0.98&0.94&0.95&0.95\\
&&20&0.97&0.98&0.98&0.97&0.97&0.97&0.95&0.96&0.97&0.89&0.90&0.90\\
&&50&0.68&0.67&0.68&0.63&0.61&0.63&0.61&0.60&0.60&0.53&0.53&0.53\\ \cline{2-15}
&\multirow{4}{*}{2}&5&1.00&1.00&1.00&1.00&1.00&1.00&0.99&1.00&1.00&0.98&0.99&0.99\\
&& 10&1.00&1.00&1.00&1.00&1.00&1.00&1.00&1.00&1.00&0.98&0.98&0.99\\
&&20&0.99&0.99&1.00&0.99&0.99&1.00&0.98&0.98&0.99&0.96&0.96&0.97\\
&& 50&0.90&0.89&0.90&0.88&0.86&0.87&0.82&0.80&0.80&0.78&0.77&0.75\\ \hline
\end{tabu}
\egroup
\caption{Fraction of calls answered in time for policies DP1, DP2 and DP3 ($K = 20$, $R = 12$, $\lambda = 0.01$).}
\label{tab::nowaiting}
\end{table}

\begin{table}[]
\small
\centering
\bgroup
\def\arraystretch{0.65}
\begin{tabu}to \textwidth {|X[c]|X[c]|X[c]|X[c]X[c]X[c]|X[c]X[c]X[c]|X[c]X[c]X[c]|X[c]X[c]X[c]|} \hline
\multirow{2}{*}{$M$}&\multirow{2}{*}{$d$}&\multirow{2}{*}{$t^*$}&\multicolumn{3}{c|}{$\mu = 0.2$}&\multicolumn{3}{c|}{$\mu = 0.1$}&\multicolumn{3}{c|}{$\mu = 0.05$}&\multicolumn{3}{c|}{$\mu = 0.02$}\\ \cline{4-15}
&&&DP1&DP4&DP5&DP1&DP4&DP5&DP1&DP4&DP5&DP1&DP4&DP5\\ \hline
\multirow{12}{*}{10}&\multirow{3}{*}{0.3}&5&0.92&0.95&0.96&0.88&0.90&0.92&0.79&0.81&0.83&0.48&0.49&0.48\\
&&10&0.80&0.91&0.90&0.74&0.85&0.87&0.61&0.73&0.77&0.31&0.41&0.42\\
&&20&0.37&0.83&0.85&0.33&0.80&0.78&0.31&0.54&0.64&0.21&0.23&0.31\\
&&50&0.24&0.64&0.54&0.23&0.52&0.52&0.22&0.43&0.41&0.21&0.22&0.21\\ \cline{2-15}
&\multirow{3}{*}{1}&5&0.97&0.97&0.98&0.94&0.94&0.96&0.86&0.87&0.88&0.63&0.59&0.59\\
&&10&0.92&0.96&0.97&0.88&0.90&0.92&0.74&0.81&0.83&0.48&0.50&0.55\\
&&20&0.60&0.93&0.91&0.59&0.83&0.86&0.44&0.63&0.78&0.31&0.30&0.43\\
&&50&0.35&0.76&0.77&0.32&0.71&0.70&0.33&0.57&0.60&0.30&0.28&0.30\\ \cline{2-15}
&\multirow{3}{*}{2}&5&0.99&0.99&0.99&0.98&0.98&0.99&0.94&0.95&0.95&0.72&0.73&0.73\\
&&10&0.97&0.98&0.98&0.94&0.96&0.97&0.89&0.91&0.92&0.65&0.64&0.66\\
&&20&0.86&0.97&0.97&0.83&0.94&0.94&0.72&0.79&0.88&0.49&0.53&0.63\\
&&50&0.54&0.88&0.87&0.50&0.84&0.87&0.45&0.71&0.72&0.43&0.34&0.43\\ \hline
\multirow{12}{*}{13}&\multirow{3}{*}{0.3}&5&0.98&0.98&0.98&0.96&0.95&0.97&0.92&0.93&0.92&0.79&0.79&0.79\\
&&10&0.94&0.97&0.97&0.93&0.94&0.94&0.88&0.90&0.91&0.69&0.73&0.73\\
&&20&0.77&0.92&0.94&0.78&0.90&0.91&0.67&0.79&0.84&0.43&0.56&0.69\\
&&50&0.31&0.83&0.77&0.31&0.77&0.77&0.30&0.69&0.70&0.28&0.41&0.56\\ \cline{2-15}
&\multirow{3}{*}{1}&5&0.99&0.99&0.99&0.98&0.99&0.99&0.95&0.96&0.97&0.86&0.86&0.86\\
&&10&0.97&0.99&0.98&0.97&0.96&0.97&0.94&0.93&0.95&0.81&0.78&0.81\\
&&20&0.91&0.96&0.96&0.88&0.95&0.96&0.80&0.89&0.92&0.62&0.67&0.80\\
&&50&0.44&0.89&0.90&0.42&0.83&0.84&0.42&0.79&0.81&0.40&0.47&0.66\\ \cline{2-15}
&\multirow{3}{*}{2}&5&1.00&0.99&0.99&0.98&0.99&0.99&0.99&0.98&0.99&0.94&0.87&0.90\\
&&10&0.99&0.98&0.99&0.99&0.98&0.99&0.97&0.96&0.97&0.91&0.85&0.91\\
&&20&0.98&0.95&0.98&0.97&0.95&0.98&0.94&0.93&0.96&0.82&0.83&0.90\\
&&50&0.64&0.93&0.91&0.64&0.93&0.91&0.58&0.90&0.93&0.53&0.61&0.82\\ \hline
\multirow{12}{*}{16}&\multirow{3}{*}{0.3}&5&0.99&0.99&0.99&0.98&0.98&0.99&0.97&0.98&0.97&0.94&0.92&0.91\\
&& 10&0.98&0.98&0.99&0.97&0.97&0.98&0.96&0.95&0.96&0.90&0.90&0.90\\
&&20&0.94&0.97&0.97&0.92&0.94&0.96&0.90&0.93&0.95&0.81&0.81&0.89\\
&&50&0.52&0.90&0.88&0.47&0.89&0.89&0.48&0.83&0.85&0.40&0.64&0.79\\ \cline{2-15}
&\multirow{3}{*}{1}&5&1.00&1.00&1.00&0.99&0.99&0.99&0.98&0.99&0.99&0.96&0.94&0.94\\
&&10&0.98&0.99&0.99&0.99&0.98&0.99&0.98&0.98&0.98&0.94&0.93&0.90\\
&&20&0.98&0.98&0.98&0.97&0.96&0.97&0.94&0.94&0.97&0.89&0.87&0.91\\
&&50&0.70&0.94&0.93&0.65&0.92&0.91&0.59&0.89&0.91&0.56&0.76&0.84\\ \cline{2-15}
&\multirow{3}{*}{2}& 5&1.00&1.00&0.99&0.99&0.99&0.99&1.00&0.98&0.98&0.98&0.97&0.95\\
&&10&1.00&0.99&0.99&1.00&0.99&0.99&1.00&0.98&0.98&0.97&0.93&0.94\\
&&20&0.99&0.98&0.98&0.99&0.96&0.99&0.99&0.95&0.97&0.95&0.90&0.94\\
&&50&0.90&0.96&0.96&0.86&0.96&0.96&0.86&0.92&0.94&0.77&0.83&0.93\\ \hline
\end{tabu}
\egroup
\caption{Fraction of calls answered in time for policies DP1, DP4 and DP5 ($K = 20$, $R = 12$, $\lambda = 0.01$).}
\label{tab::waiting}
\end{table}

\begin{table}[]
\centering
\bgroup
\def\arraystretch{0.7}
\begin{tabu}to \textwidth { | X[c] |X[c]|X[c]|X[c]X[c]X[c]X[c]X[c]| } 
\hline
$\mu$& $t^*$ &$d$&RP1&RP2&RP3&RP4&RP5\\ \hline
\multirow{9}{*}{0.2}&\multirow{3}{*}{5}&0.3&0.93&0.69&0.71&0.89&0.98\\
&&1&0.94&0.79&0.81&0.95&0.99\\
&&2&0.97&0.87&0.89&0.99&0.99\\ \cline{2-8}
&\multirow{3}{*}{10}&0.3&0.86&0.58&0.57&0.72&0.95\\
&&1&0.91&0.74&0.75&0.85&0.97\\
&&2&0.94&0.82&0.87&0.94&0.99\\ \cline{2-8}
&\multirow{3}{*}{20}&0.3&0.77&0.43&0.43&0.48&0.92\\
&&1&0.83&0.61&0.62&0.64&0.96\\
&&2&0.92&0.74&0.79&0.79&0.97\\ \cline{1-8}
\multirow{9}{*}{0.1}&\multirow{3}{*}{5}&0.3&0.86&0.70&0.73&0.91&0.96\\
&&1&0.87&0.84&0.85&0.96&0.98\\
&&2&0.94&0.92&0.95&0.98&0.99\\ \cline{2-8}
&\multirow{3}{*}{10}&0.3&0.83&0.57&0.56&0.69&0.94\\
&&1&0.84&0.74&0.74&0.81&0.96\\
&&2&0.88&0.81&0.85&0.94&0.99\\ \cline{2-8}
&\multirow{3}{*}{20}&0.3&0.70&0.32&0.31&0.47&0.91\\
&&1&0.75&0.58&0.58&0.60&0.95\\
&&2&0.83&0.69&0.73&0.75&0.97\\ \cline{1-8}
\multirow{9}{*}{0.05}&\multirow{3}{*}{5}&0.3&0.76&0.70&0.70&0.85&0.95\\
&&1&0.80&0.83&0.86&0.95&0.97\\
&&2&0.87&0.92&0.96&0.98&0.99\\ \cline{2-8}
&\multirow{3}{*}{10}&0.3&0.66&0.46&0.47&0.68&0.91\\
&&1&0.71&0.64&0.64&0.81&0.95\\
&&2&0.77&0.82&0.86&0.93&0.98\\ \cline{2-8}
&\multirow{3}{*}{20}&0.3&0.52&0.25&0.25&0.43&0.82\\
&&1&0.56&0.37&0.38&0.58&0.91\\
&&2&0.64&0.49&0.55&0.75&0.96\\ \hline
\end{tabu}
\egroup
\caption{Fraction of calls answered in time for policies RP1-RP5 ($K = 20$, $R = 10$, $M = 13$, $\lambda = 0.01$).}
\label{tab::combpolicy}
\end{table}


\end{document}